\documentclass[12pt]{article}

\usepackage[english]{babel}

\usepackage[a4paper]{geometry}
\usepackage{amssymb,amsmath,amsthm}
\usepackage{graphicx,cite,color}
\usepackage{longtable,pdflscape,booktabs,caption,multicol}
\usepackage[colorlinks=true,citecolor=black,linkcolor=black,urlcolor=blue]{hyperref}
\usepackage{enumitem}
\usepackage{mathtools}

\theoremstyle{plain}
\newtheorem{theorem}{Theorem}
\newtheorem{lemma}[theorem]{Lemma}

\newtheorem{proposition}[theorem]{Proposition}

\theoremstyle{definition}

\theoremstyle{remark}
\newtheorem*{remark}{Remark}

\newcommand{\arxiv}[1]{arXiv: \href{https://arxiv.org/abs/#1}{\texttt{#1}}}
\newcommand{\doi}[1]{DOI: \href{https://doi.org/#1}{\texttt{#1}}}

\def\Z{\mathbb{Z}}

\usepackage{float}
\usepackage{tikz}
\usepackage{subfig}
\usetikzlibrary{automata}
\usetikzlibrary{decorations.pathmorphing}
\tikzset{snake it/.style={decorate, decoration=snake}}
\usetikzlibrary{arrows}
\usetikzlibrary{arrows.meta}
\usetikzlibrary{chains}
\usetikzlibrary{positioning}
\usetikzlibrary{automata,positioning,calc}
\usetikzlibrary{decorations}
\usetikzlibrary{decorations.shapes}
\usetikzlibrary{decorations.markings}

\tikzset{
    edge/.style={-{Latex[scale=1.7]}},
}

\usepackage{blkarray, bigstrut}
\usepackage{bm}

\usepackage{listings}
\definecolor{mauve}{rgb}{0.58,0,0.82}
\lstdefinestyle{pitonche} {
    language = Python,
    basicstyle = footnotesizettfamily,
    showspaces = false,
    showstringspaces = false,
    breakautoindent = true,
    flexiblecolumns = true,
    keepspaces = true,
    stepnumber = 1,
    xleftmargin = 0pt
}
\def\dj{d\kern-0.4em\char"16\kern-0.1em}
\def\Dj{\hbox{\raise0.3ex\hbox{-}\kern-0.4em  D}}

\usepackage{authblk}

\title{Classification of cubic tricirculant nut graphs}

\author[1,2]{Ivan Damnjanović\thanks{The corresponding author.}}
\author[3,4,5]{Nino Bašić}
\author[3,4,5]{Tomaž Pisanski}
\author[5,6]{Arjana Žitnik}
\affil[1]{Faculty of Electronic Engineering, University of Niš, Niš, Serbia}
\affil[2]{Diffine LLC, San Diego, California, USA}
\affil[3]{FAMNIT, University of Primorska, Koper, Slovenia}
\affil[4]{IAM, University of Primorska, Koper, Slovenia}
\affil[5]{Institute of Mathematics, Physics and Mechanics, Ljubljana, Slovenia}
\affil[6]{Faculty of Mathematics and Physics, University of Ljubljana, Ljubljana, Slovenia}

\date{}

\begin{document}

\maketitle

\begin{abstract}
A nut graph is a simple graph whose adjacency matrix has the eigenvalue zero
with multiplicity one such that its corresponding eigenvector has no zero entries.
It is known that there exist no cubic circulant nut graphs. 
A bicirculant (resp.\ tricirculant) graph is defined as a graph that admits a cyclic group
of automorphisms having two (resp.\ three) orbits of vertices of equal size.
We show that there   exist no cubic bicirculant nut graphs and we provide a full classification of  cubic tricirculant nut graphs.
\end{abstract}


\section{Introduction}


A \emph{nut graph} is a graph of nullity one whose null space is spanned by a full vector, i.e., a vector with no zero entries. Nut graphs were introduced and studied by Sciriha and Gutman \cite{Sciriha1997, Sciriha1998_A, Sciriha1998_B, ScGu1998, Sciriha1999}. Some of their properties were further investigated in \cite{Sciriha2007, Sciriha2008, FoGaGoPiSc2020, GaPiSc2023}. Moreover, the chemical justification for studying such graphs can be found in \cite{ScFo2007, ScFo2008, FoPiToBoSc2014, CoFoGo2018, FoPiBa2021} and many other results concerning them 
are to be found in the monograph \cite{ScFa2021}.

The problem of the existence of regular nut graphs on a given number of vertices and a given degree  was first considered in \cite{FoGaGoPiSc2020}, where it was solved for degrees up to 11. In \cite{BaKnSk2022} the degree was extended to 12 using circulant graphs in a clever way.
%
%
In a series of papers \cite{DaSt2022, Damnjanovic2022_A, Damnjanovic2022_B}, Damnjanović and Stevanović completely resolved the problem of existence of circulant nut graphs of given degree and order by establishing the following result.

\begin{theorem}[\hspace{1sp}{\cite[Theorem 5]{Damnjanovic2022_B}}]\label{damnjanovic_th}
Let $d$ and $n$ be nonnegative integers. Then there exists a $d$-regular
circulant nut graph of order $n$ if and only if $d > 0$, $4 \mid d$, $2 \mid n$, 
together with $n \ge d + 4$ if $d \equiv_8 4$, and $n \ge d + 6$ if $8 \mid d$, 
as well as $(n, d) \neq (16, 8)$.
\end{theorem}


From Theorem~\ref{damnjanovic_th} it  follows that there are no cubic circulant nut graphs and that there exists a quartic circulant nut graph of every even order $n \ge 8$. Moreover, the full characterization of quartic circulant nut graphs can be found in \cite{Damnjanovic2022_C}.

One can generalize the concept of a circulant graph by taking into consideration the so-called \emph{bircirculant} (resp.\ \emph{tricirculant}) graphs, which are the graphs that  admit  a cyclic group of automorphisms having two (resp.\ three) orbits of vertices of equal size  \cite{KoKuMaWi2012, Pi2007}. 
In this paper we show that there exist no cubic bircirculant nut graphs. 
While  there are no cubic circulant and bicirculant graphs, there exist cubic tricirculant graphs. 
We give  a complete classification of these graphs. 
Below is a more detailed description of our classification.
We also note that the quartic bircirculant nut graphs are currently being investigated \cite{PiPoZi2022}.

For convenience, we will suppose that the vertex set of each tricirculant graph of order $3n$ is given via $X \cup Y \cup Z$, where
\[
    X = \{ x_0, x_1, \ldots, x_{n-1} \}, \qquad Y = \{ y_0, y_1, \ldots, y_{n-1} \}, \qquad Z = \{ z_0, z_1, \ldots, z_{n-1} \} ,
\]
and that there exists an automorphism which maps $x_j$ to $x_{j+1}$, $y_j$ to $y_{j+1}$ and $z_j$ to $z_{j+1}$, for each $j \in \mathbb{Z}_n$. If $G$ is an arbitrary cubic tricirculant graph of order $3n$, it is then clear that $n$ is necessarily even. Moreover, as observed by Potočnik and Toledo \cite[Theorem 2.4]{PoTol2020}, it can be shown that if $G$ is connected, then it must be isomorphic to a graph that belongs to at least one of the next four families:
\begin{enumerate}[label={\bf (\roman*)}]
\item $T_1(n, a, b),\, 0 \le a < b < n$;
\item $T_2(n, a, b),\, 0 < a < n,\, 0 < b < \frac{n}{2}$;
\item $T_3(n, a),\, 0 \leq a < n$;
\item $T_4(n, a, b),\, 0 < a \leq b < \frac{n}{2}$;
\end{enumerate}
whose elements can be defined via their edge sets as follows:
\begin{align*}
E(T_1(n,a,b)) &= \{ x_j y_{j+a},\, x_j y_{j+b},\, x_j z_j,\, y_j z_j,\, z_j z_{j+\frac{n}{2}} \mid j \in \Z_n \},\\
E(T_2(n,a,b)) &= \{ x_j x_{j+b},\, x_j z_j,\, y_j y_{j+\frac{n}{2}},\,  y_j z_j,\, y_j z_{j-a} \mid j \in \Z_n \},\\
E(T_3(n,a)) & = \{ x_j x_{j+\frac{n}{2}},\, x_j y_{j+a},\,  x_j z_j,\, y_j y_{j+\frac{n}{2}},\, y_j z_j,\, z_j z_{j+\frac{n}{2}} \mid j \in \Z_n \}, \\
E(T_4(n, a, b)) & = \{ x_j x_{j+a},\, x_j z_j,\, y_j y_{j+b},\, y_j z_j,\, z_j z_{j+\frac{n}{2}} \mid j \in \Z_n \}.
\end{align*}

For brevity, if a cubic tricirculant graph is isomorphic to at least one element from the $T_j$ family, we will then say that such a graph is of \emph{type $j$}, for each $j = 1, 2, 3, 4$. 
Note that in principle, a cubic tricirculant could be of more than one type. However,
every vertex-transitive cubic tricirculant is of type $j$ for exactly one $j \in \{1,2,3,4\}$, with
the exception of the triangular prism which is both of type~1 as well as of type~3 \cite{PoTol2020}.
The above definitions can be concisely visualized in the form of Figure \ref{voltagegraphs}. Here, we may note that this figure actually formally represents the \emph{voltage graphs} corresponding to the cubic tricirculants \cite{PoTol2020}.

\begin{figure}[htb]
\centering
\subfloat[Type $1$]{ 
    \scalebox{0.65}{\begin{tikzpicture}
        \node[state, minimum size=0.75cm, thick] (1) at (1.25, -2.17) {$X$};
        \node[state, minimum size=0.75cm, thick] (2) at (2.5, 0) {$Z$};
        \node[state, minimum size=0.75cm, thick] (3) at (3.75, -2.17) {$Y$};
        \node[state, minimum size=0.75cm, color=white] (4) at (2.5, 2) {$D$};
        \node[state, minimum size=0.75cm, color=white] (6) at (4.75, -4) {$F$};

        \path[thick] (1) edge node[pos=0.4, above, xshift=-0.2cm] {$0$} (2);
        \path[thick] (2) edge node[pos=0.6, above, xshift=0.2cm] {$0$} (3);
        \draw[edge] (1) to[bend left=25] node[pos=0.5, above] {$a$} (3);
        \draw[edge] (1) to[bend right=25] node[pos=0.5, below] {$b$} (3);
        \path[thick] (2) edge node[pos=1, right] {$\dfrac{n}{2}$} (4);
    \end{tikzpicture}
}}
\hspace{-0.1cm}
\subfloat[Type $2$]{ 
\scalebox{0.65}{    \begin{tikzpicture}
        \node[state, minimum size=0.75cm, thick] (1) at (1.25, -2.17) {$X$};
        \node[state, minimum size=0.75cm, thick] (2) at (2.5, 0) {$Z$};
        \node[state, minimum size=0.75cm, thick] (3) at (3.75, -2.17) {$Y$};
        \node[state, minimum size=0.75cm, color=white] (5) at (0.25, -4) {$E$};
        \node[state, minimum size=0.75cm, color=white] (6) at (4.75, -4) {$F$};

        \draw[edge] (1) to[out=135, in=210, looseness=6] node[pos=0.2, above] {$b$} (1);
        \path[thick] (3) edge node[pos=0.5, yshift=0.2cm, right] {$\dfrac{n}{2}$} (6);
        \path[thick] (1) edge node[pos=0.5, above, xshift=-0.2cm] {$0$} (2);
        \draw[edge] (2) to[bend left=25] node[pos=0.5, above, xshift=0.2cm] {$0$} (3);
        \draw[edge] (2) to[bend right=25] node[pos=0.5, below, xshift=-0.2cm] {$a$} (3);
    \end{tikzpicture}
}}
\hspace{-0.1cm}
\subfloat[Type $3$]{
\scalebox{0.65}{    \begin{tikzpicture}
        \node[state, minimum size=0.75cm, thick] (1) at (1.25, -2.17) {$X$};
        \node[state, minimum size=0.75cm, thick] (2) at (2.5, 0) {$Z$};
        \node[state, minimum size=0.75cm, thick] (3) at (3.75, -2.17) {$Y$};
        \node[state, minimum size=0.75cm, color=white] (4) at (2.5, 2) {$D$};
        \node[state, minimum size=0.75cm, color=white] (5) at (0.25, -4) {$E$};
        \node[state, minimum size=0.75cm, color=white] (6) at (4.75, -4) {$F$};

        \path[thick] (1) edge node[pos=0.4, above, xshift=-0.2cm] {$0$} (2);
        \path[thick] (2) edge node[pos=0.6, above, xshift=0.2cm] {$0$} (3);
        \draw[edge] (1) to node[pos=0.5, above] {$a$} (3);
        \path[thick] (2) edge node[pos=1, right] {$\dfrac{n}{2}$} (4);
        \path[thick] (1) edge node[pos=0.5, yshift=0.2cm, left] {$\dfrac{n}{2}$} (5);
        \path[thick] (3) edge node[pos=0.5, yshift=0.2cm, right] {$\dfrac{n}{2}$} (6);

    \end{tikzpicture}
}}
\hspace{-0.1cm}
\subfloat[Type $4$]{
\scalebox{0.65}{    \begin{tikzpicture}
        \node[state, minimum size=0.75cm, thick] (1) at (1.25, -2.17) {$X$};
        \node[state, minimum size=0.75cm, thick] (2) at (2.5, 0) {$Z$};
        \node[state, minimum size=0.75cm, thick] (3) at (3.75, -2.17) {$Y$};
        \node[state, minimum size=0.75cm, color=white] (4) at (2.5, 2) {$D$};
        \node[state, minimum size=0.75cm, color=white] (5) at (0.25, -4) {$E$};
        \node[state, minimum size=0.75cm, color=white] (6) at (4.75, -4) {$F$};

        \path[thick] (1) edge node[pos=0.4, above, xshift=-0.2cm] {$0$} (2);
        \path[thick] (2) edge node[pos=0.6, above, xshift=0.2cm] {$0$} (3);
        \draw[edge] (1) to[out=135, in=210, looseness=6] node[pos=0.2, above] {$a$} (1);
        \draw[edge] (3) to[out=45, in=-30, looseness=6] node[pos=0.2, above] {$b$} (3);
        \path[thick] (2) edge node[pos=1, right] {$\dfrac{n}{2}$} (4);
    \end{tikzpicture}
}}
\caption{Voltage graphs for the cubic tricirculant graphs
of types $1$, $2$, $3$ and $4$.} 
\label{voltagegraphs}
\end{figure}

Finally, for each $x \in \mathbb{Z}, \, x \neq 0$, let $v_2(x)$ denote the power of two in the prime factorization of $|x|$, i.e., the unique $\beta \in \mathbb{N}_0$ such that $2^\beta \mid x$, but $2^{\beta+1} \nmid x$. We now state our main theorem.

\begin{theorem}[Tricirculant cubic nut graph classification]\label{main_theorem}
A tricirculant cubic graph is a nut graph if and only if it is representable as a $T_1(n, a, b), \, 2 \mid n, \, 0 \le a < b < n$ such that
\begin{enumerate}[label=\textbf{(\roman*)}]
    \item $\gcd\left( \frac{n}{2}, a \right) = \gcd\left( \frac{n}{2}, b \right) = 1$;
    \item $a \not\equiv_2 \frac{n}{2}$ and $b \not\equiv_2 \frac{n}{2}$;
    \item $v_2(b - a) \ge v_2(n)$;
\end{enumerate}
or as a $T_4(n, a, b), \, 2 \mid n, \, 0 < a \le b < \frac{n}{2}$ where
\begin{enumerate}[label=\textbf{(\roman*)}]\setcounter{enumi}{3}
    \item $\gcd\left( \frac{n}{2}, a, b \right) = 1$;
    \item if $4 \nmid n$, then at least one of $a, b$ is even;
    \item if $4 \mid n$, then $a$ and $b$ are of different parities;
    \item if $10 \mid n$, then at least one of $a, b, a-b, a+b$ is divisible by five.
\end{enumerate}
\end{theorem}

\begin{remark}
    Condition \textbf{(iv)} from Theorem \ref{main_theorem} is actually equivalent to $T_4(n, a, b)$ being a connected graph (see, for example, \cite[Theorem 2.4]{PoTol2020}). Moreover, condition \textbf{(i)} also implies that $T_1(n, a, b)$ is surely connected.
\end{remark}

The rest of the paper will focus on providing the full proof of Theorem \ref{main_theorem}. Bearing this in mind, its structure shall be organized as follows. Section \ref{sc_preliminaries} will serve to preview certain theoretical results to be used throughout the remaining sections. Afterwards, we will use Section~\ref{bicirc_sc} to show that no cubic bicirculant graph can be a nut graph and in Section~\ref{sc_two_families} we will take into consideration an arbitrary tricirculant cubic nut graph and demonstrate that it is necessarily of type $1$ or of type $4$. Afterwards, Section \ref{sc_type_b} will be used to obtain the precise conditions that a type $1$ graph should satisfy in order to be a nut graph. Finally, we shall determine all the type $4$ nut graphs in Section \ref{sc_type_a}, thereby completing the proof of Theorem~\ref{main_theorem}.

\section{Preliminaries}\label{sc_preliminaries}

In this section we review some known theoretical results from various fields of mathematics which will be used later on throughout the remaining sections. First of all, recall that $\eta(G)$ denotes the \emph{nullity} of $G$, i.e., the nullity of the adjacency matrix of the graph $G$. Similarly, we will use $\mathcal{N}(A_G)$ to denote the null space of $G$. We also note that the eigenvectors corresponding to the eigenvalue zero are also called the \emph{kernel   vectors} and may be characterized via the following lemma.

\begin{lemma}\label{local_condition_lemma}
Let $G$ be a  graph and let $u \in \mathbb{R}^{V_G}$ be a nonzero vector. Then $u$ is an eigenvector  corresponding to the eigenvalue zero if and only if for every vertex $x$ of $G$ the following holds:
\[
    \sum_{y \sim x} u(y)=0.
\]
\end{lemma}
\begin{proof}
Multiply the row of the adjacency matrix of $G$, corresponding to the vertex $x$,
by the eigenvector $u$, and the result is obtained.
\end{proof} 

Here, we shall call the above condition that a kernel vector must satisfy the \emph{local condition}. We now recall some basic properties of nut graphs, see \cite{ScGu1998}.

\begin{lemma}\label{lemma_conn}
Every nut graph is connected.
\end{lemma}

\begin{lemma}  \label{lemma_bip}
A bipartite graph is not a nut graph.
\end{lemma}

Furthermore, we will present the following lemma which forms a connection between the null space vectors of a nut graph and the vertices of a particular orbit of a given graph automorphism. For the proof, see, for example, \cite[p.\ 135]{Cvetkovic1995}.
\begin{lemma}\label{nut_orbit_lemma}
Let $G$ be a nut graph and let $\pi \in \mathrm{Aut}(G)$ be its automorphism. 
If $u \in \mathcal{N}(A_G)$ and $X = \{ x_0, x_1, \ldots, x_{k-1} \} \subseteq V(G)$ 
represents an orbit of $\pi$ such that
    \[
        \pi(x_j) = x_{j+1} \qquad (j = 0, \dots, k-1),
    \]
where the addition is done modulo $k$, then we have that $u$ is constant on $X$,
or the orbit size $k$ is even and
    \[
        u(x_j) = (-1)^j \, u(x_0) \qquad  (j = 0, \dots, k-1).
    \]
\end{lemma}

Circulant matrices will also play an important role in our proofs. We note that a circulant matrix $C \in \mathbb{R}^{n \times n}$ is any matrix bearing the form
\[
    C= \begin{bmatrix}
        c_0 & c_1 & c_2 & \cdots & c_{n-1}\\
        c_{n-1} & c_0 & c_1 & \cdots & c_{n-2}\\
        c_{n-2} & c_{n-1} & c_0 & \cdots & c_{n-3}\\
        \vdots & \vdots & \vdots & \ddots & \vdots\\
        c_1 & c_2 & c_3 & \dots & c_0
    \end{bmatrix} .
\]
It is well known from elementary linear algebra theory (see, for example, \cite[Section~3.1]{Gray}) that the eigenvalues of a circulant matrix can be evaluated by applying the expression
\begin{equation}\label{circulant_eigenvalues}
    c_0 + c_1 \zeta + c_2 \zeta^2 + \cdots + c_{n-1} \zeta^{n-1},
\end{equation}
as $\zeta$ ranges through the $n$-th roots of unity.

Finally, we recall the cyclotomic polynomials and provide a theorem on their divisibility which shall play a key role in Section \ref{sc_type_a}. The cyclotomic polynomial $\Phi_f(x)$ can be defined for each $f \in \mathbb{N}$ via
\[
    \Phi_f(x) = \prod_{\xi} (x - \xi) ,
\]
where $\xi$ ranges over the primitive $f$-th roots of unity. It is known that these polynomials have integer coefficients and that they are all irreducible in $\mathbb{Q}[x]$ (see, for example, \cite{Cyclotomic}). For this reason, an arbitrary polynomial in $\mathbb{Q}[x]$ has a primitive $f$-th root of unity among its roots if and only if it is divisible by $\Phi_f(x)$. Besides that, it is worth pointing out that
\begin{equation}\label{cyclotomic_reduction}
    \Phi_f(x) = \Phi_{f / p^{k-1}}(x^{p^{k-1}})
\end{equation}
holds whenever $p^k \mid f$ for a given prime number $p$ and some $k \in \mathbb{N}, \, k \ge 2$ (see, for example, \cite[p.\ 160]{Nagell1951}). We end the section by disclosing the following theorem on the divisibility of lacunary polynomials by cyclotomic polynomials.
\begin{theorem}[Filaseta and Schinzel \cite{FiSc2003}]\label{filaseta}
Let $P(x) \in \mathbb{Z}[x]$ have $N$ nonzero terms and let $\Phi_f(x) \mid P(x)$.
Suppose that $p_1, p_2, \dots, p_k$ are distinct primes such that
\[
    \sum_{j=1}^k (p_j-2) > N-2 .
\]
Let $e_j$ be the largest exponent such that $p_j^{e_j} \mid f$. Then for at least one $j$, $1 \le j \le k$, we have that $\Phi_{f'}(x)\mid P(x)$, where $f' = \dfrac{f}{p_j^{e_j}}$.
\end{theorem}

\section{Nonexistence of cubic bicirculant nut graphs}\label{bicirc_sc}

In this section we show that a cubic bicirculant graph cannot be a nut graph. As demonstrated by Pisanski \cite{Pi2007}, it is not difficult to establish that each connected cubic bicirculant graph of order $2n$ must be isomorphic to a graph from at least one of the following three families:
\begin{enumerate}[label={\bf (\roman*)}]
\item $B_1(n, a, b),\, 0 < a < b < n$;
\item $B_2(n, a),\, 0 < a < n$;
\item $B_3(n, a, b),\, 0 < a \le b < \frac{n}{2}$;
\end{enumerate}
whose elements have the vertex set $\{x_0, \ldots, x_{n-1}, y_0, \ldots, y_{n-1} \}$ and the edge sets as given below:
\begin{align*}
E(B_1(n, a, b)) &= \{ x_j y_j,\, x_j y_{j+a},\, x_j y_{j+b} \mid j \in \mathbb{Z}_n \},\\
E(B_2(n, a)) &= \{ x_j x_{j+\frac{n}{2}},\, x_j y_j,\, x_j y_{j+a}, \, y_j y_{j+\frac{n}{2}} \mid j \in \mathbb{Z}_n \},\\
E(B_3(n, a, b)) &= \{ x_j y_j,\, x_j x_{j+a},\, y_j y_{j+b} \mid j \in \Z_n \}.
\end{align*}

\begin{figure}[htb]
\centering
\subfloat[Type 1]{ 
    \scalebox{0.65}{\begin{tikzpicture}
        \node[state, minimum size=0.75cm, thick] (1) at (0, 0) {$X$};
        \node[state, minimum size=0.75cm, thick] (3) at (3, 0) {$Y$};

        \draw[edge] (1) to[bend left=45] node[pos=0.5, above] {$a$} (3);
        \draw[edge] (1) to[bend right=45] node[pos=0.5, below] {$b$} (3);
        \path[thick] (1) edge node[above] {$0$} (3);
    \end{tikzpicture}
}}
\hspace{-0.1cm}
\subfloat[Type 2]{ 
\scalebox{0.65}{    \begin{tikzpicture}
        \node[state, minimum size=0.75cm, thick] (1) at (0, 0) {$X$};
        \node[state, minimum size=0.75cm, thick] (3) at (3, 0) {$Y$};
        \node[state, minimum size=0.75cm, thick, color=white] (5) at (-1.5, 0) {$E$};
        \node[state, minimum size=0.75cm, thick, color=white] (6) at (4.5, 0) {$F$};

        \draw[edge] (1) to[bend left=25] node[pos=0.5, above] {$a$} (3);
        \path[thick] (1) edge[bend right=25] node[below] {$0$} (3);
        \path[thick] (1) edge node[pos=0.5, yshift=0.2cm, above] {$\dfrac{n}{2}$} (5);
        \path[thick] (3) edge node[pos=0.5, yshift=0.2cm, above] {$\dfrac{n}{2}$} (6);
    \end{tikzpicture}
}}
\hspace{-0.1cm}
\subfloat[Type 3]{
\scalebox{0.65}{    \begin{tikzpicture}
        \node[state, minimum size=0.75cm, thick] (1) at (0, 0) {$X$};
        \node[state, minimum size=0.75cm, thick] (3) at (3, 0) {$Y$};

        \draw[edge] (1) to[out=135, in=210, looseness=6] node[pos=0.2, above] {$a$} (1);
        \draw[edge] (3) to[out=45, in=330, looseness=6] node[pos=0.2, above] {$b$} (3);
        \path[thick] (1) edge node[above] {$0$} (3);
    \end{tikzpicture}
}}
\caption{Voltage graphs for the cubic bicirculant graphs of types 1, 2 and 3.} 
\label{voltagegraphs2}
\end{figure}
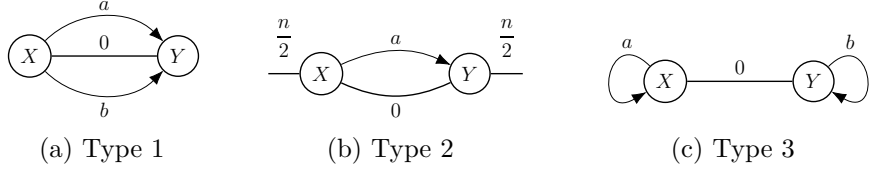

To make matters simple, we will say that a cubic bicirculant graph that is isomorphic to at least one element from the $B_j$ family is of \emph{type $j$}, for each $j = 1, 2, 3$. This terminology can now be concisely visualized in the form of Figure \ref{voltagegraphs2}. In the following proposition, we will demonstrate that truly none of these graphs can be nut graphs.

\begin{proposition}
    There does not exist a cubic bicirculant nut graph.
\end{proposition}
\begin{proof}
    Due to Lemma \ref{lemma_conn}, it is sufficient to restrict ourselves to connected graphs. Let $G$ be a connected cubic bircirculant graph. If $G$ is of type $1$, then this graph is necessarily bipartite, hence it cannot be a nut graph, by virtue of Lemma \ref{lemma_bip}. If we suppose that $G$ is of type~$2$, then applying the local conditions on the vector $u \in \mathcal{N}(A_G)$ promptly gives us
    \begin{align}
    \label{auxx_1} u(x_{j + \frac{n}{2}}) + u(y_j) + u(y_{j+a}) &= 0 \qquad (j \in \Z_n) ,\\
    \label{auxx_2} u(x_j) + u(x_{j-a}) + u(y_{j + \frac{n}{2}}) &= 0 \qquad (j \in \Z_n) .
    \end{align}
    By virtue of Lemma \ref{nut_orbit_lemma}, we see that $u(y_j) = u(y_{j+a})$ or $u(y_j) = -u(y_{j+a})$, which means that Equation \eqref{auxx_1} leads us to $u(x_{j+\frac{n}{2}}) = 0$ or $u(x_{j+\frac{n}{2}}) = -2u(y_j)$. The former would immediately imply that $G$ is not a nut graph, hence we may assume that the latter holds. Analogously, from Equation \eqref{auxx_2} we can conclude that $u(y_{j+\frac{n}{2}}) = 0$ or $u(y_{j+\frac{n}{2}}) = -2u(x_j)$. From $u(y_{j+\frac{n}{2}}) = 0$ we again obtain that $G$ is not a nut graph, while the latter quickly yields $u(x_j) = 4 u(x_j)$. From here, it promptly follows that $u(x_j) = 0$.

    Finally, if we take $G$ to be of type $3$, then the local conditions on the vector $u \in \mathcal{N}(A_G)$ dictate
    \begin{align}
    \label{auxx_3} u(x_{j+a}) + u(x_{j-a}) + u(y_j) &= 0 \qquad (j \in \Z_n) ,\\
    \label{auxx_4} u(x_j) + u(y_{j+b}) + u(y_{j - b}) &= 0 \qquad (j \in \Z_n) .
    \end{align}
    Lemma \ref{nut_orbit_lemma} guarantees that $u(x_{j-a}) = u(x_{j+a})$ and $u(y_{j-b}) = u(y_{j+b})$ must hold, which means that Equations \eqref{auxx_3} and \eqref{auxx_4} directly imply $u(y_j) = -2u(x_{j+a})$ and $u(x_j) = -2u(y_{j+b})$, respectively. However, this quickly gives us $u(x_j) = 4 u(x_j)$ or $u(x_j) = -4 u(x_j)$, hence $u(x_j) = 0$. 
\end{proof}

\section[Nonexistence of tricirculant nut graphs of types 2 and 3]{Nonexistence of tricirculant nut graphs\\ of types 2 and 3}\label{sc_two_families}

In this section, we will give a brief demonstration that if a tricirculant cubic graph is a nut graph, then it must necessarily belong to the family comprising the type~1 graphs or the family consisting of the type~4 graphs. 

\begin{proposition}
A tricirculant cubic graph of type $2$ or type $3$ is not a nut graph.
\end{proposition}
\begin{proof}
Let $G$ be a tricirculant cubic nut graph of type $2$, say $G$ is isomorphic to $T_2(n, a, b)$ for appropriate parameters $n,a,b$, and let $u \in \mathcal{N}(A_G)$. 
Thus, we may assume that $G$ has the vertex set $\{x_0,\dots,x_{n-1},y_0,\dots,y_{n-1},z_0,\dots,z_{n-1}\}$
and that $x_j \sim  x_{j+b},x_{j-b}, z_j$, $y_j \sim y_{j+\frac{n}{2}}, z_j, z_{j-a}$ and
$z_j \sim x_j, y_j, y_{j+a}$ for $j \in \Z_n$.
Taking everything into consideration, we see that the local conditions for the vector $u$ must bear the form
\begin{align}
\label{aux_3} u(x_{j+b}) + u(x_{j-b}) + u(z_j) &= 0 \qquad (j \in \Z_n),\\
\label{aux_2} u(y_{j + \frac{n}{2}}) + u(z_j) + u(z_{j-a}) &= 0 \qquad (j \in \Z_n),\\
\label{aux_4} u(y_j) + u(y_{j+a}) + u(x_j) &= 0 \qquad (j \in \Z_n),
\end{align}
where Equations \eqref{aux_3}, \eqref{aux_2}, \eqref{aux_4} represent the local conditions at vertices $x_j$, $y_j$ and $z_j$, respectively.
    
    By applying Lemma \ref{nut_orbit_lemma} multiple times, it is not difficult to reach a contradiction. For starters, given the fact that $j+b \equiv_2 j-b$, Equation\ (\ref{aux_3}) immediately tells us that $u(z_j) \in \{ 2 u(x_j), -2 u(x_j) \}$. Furthermore, if $u(z_j) = -u(z_{j-a})$ (resp.\ $u(y_j) = -u(y_{j+a})$), then Equation\ (\ref{aux_2}) (resp.\ (\ref{aux_4})) yields $u(y_j) = 0$ (resp.\ $u(x_j) = 0$), which leads us to $u = \bm{0}$. On the other hand, if $u(z_j) = u(z_{j-a})$ and $u(y_j) = u(y_{j+a})$ both hold, then Equation\ (\ref{aux_2}) implies $u(y_j) \in \{ 4u(x_j), -4u(x_j) \}$, while Equation~(\ref{aux_4}) subsequently gives us $u(x_j) \in \{8u(x_j), -8u(x_j) \}$. For this reason, $u(x_j) = 0$ must hold, which promptly implies $u = \bm{0}$ once again. Thus, in any case, the zero vector is certainly the only element of $\mathcal{N}(A_G)$, hence $\eta(G) = 0$, which leads to a contradiction.

Now let $G$ be a tricirculant cubic nut graph of type 3, say $G$ is isomorphic to $T_3(n,a)$ for appropriate parameters $n,a$. Thus, we may assume that $G$ has the vertex set $\{x_0,\dots,x_{n-1},y_0,\dots,y_{n-1},z_0,\dots,z_{n-1}\}$ and that $x_j \sim  x_{j+\frac{n}{2}}, y_{j+a}, z_j$, $y_j \sim x_{j-a},y_{j+\frac{n}{2}}, z_{j}$ and $z_j \sim x_j, y_j, z_{j+\frac{n}{2}}$ for $j \in \Z_n$.
Now let us define two vectors $u, v \in \mathbb{R}^{V_G}$ such that
    \begin{alignat*}{5}
        u(x_j) &= 1, \quad & u(y_j) &= -1, \quad & u(z_j) &= 0 \qquad && (j \in \Z_n),\\
        v(x_j) &= 1, \quad & v(y_j) &= 0, \quad & v(z_j) &= -1  \qquad &&  (j \in \Z_n).
    \end{alignat*}
Then we quickly see that these two vectors are linearly independent, while both of them belong to $\mathcal{N}(A_G)$. Hence, $\eta(G) \ge 2$, which means that $G$ cannot be a nut graph.
\end{proof}

\section{Tricirculant graphs of type 1}\label{sc_type_b}

In the previous section we have demonstrated that every tricirculant cubic nut graph must be of type 1 or type 4. Our next step in proving Theorem~\ref{main_theorem} will be to precisely determine all the nut graphs among the tricirculant cubic graphs of   type~1. In order to achieve this, we will use the given section to show that the following theorem holds.

\begin{theorem}\label{type_b_th}
An arbitrary graph representable as $T_1(n, a, b)$, where $n$ is even and $0 \le a, b < n, \, a \neq b$, is a nut graph if and only if the following conditions hold:
    \begin{enumerate}[label={\bf (\roman*)}]
        \item $\gcd\left( \frac{n}{2}, a \right) = \gcd\left( \frac{n}{2}, b \right) = 1$;
        \item $a \not\equiv_2 \frac{n}{2}$ and $b \not\equiv_2 \frac{n}{2}$;
        \item $v_2(b - a) \ge v_2(n)$.
    \end{enumerate}
\end{theorem}

Before we give the proof of Theorem \ref{type_b_th} itself, we will need one auxiliary claim that connects the nut property of $T_1(n, a, b)$ to the root properties of a concrete polynomial. The next lemma demonstrates the aforementioned observation.

\begin{lemma}\label{type_b_polynomial_lemma}
    A graph representable as $T_1(n, a, b)$ is a nut graph if and only if the $\mathbb{Z}[x]$ polynomial
    \begin{equation}\label{aux_12}
         x^{2a + b} + x^{a + 2b} + x^a + x^b - x^{\frac{n}{2} + 2a} - x^{\frac{n}{2} + 2b} - 2x^{\frac{n}{2} + a + b}
    \end{equation}
    has no $n$-th roots of unity among its roots, besides $1$.
\end{lemma}
\begin{proof}
Let $G$ be a given graph which is representable as $T_1(n, a, b)$. Thus, we may assume that $G$ has the vertex set 
$\{x_0,\dots,x_{n-1},y_0,\dots,y_{n-1},z_0,\dots,z_{n-1}\}$
and that $x_j \sim y_{j+a}, y_{j+b}, z_j$, $y_j \sim x_{j-a}, x_{j-b}, z_{j}$ and $z_j \sim x_j,y_j, z_{j+\frac{n}{2}}$ for $j \in \Z_n$.

From the local conditions for $x_j,y_j,z_j$, respectively, it is clear that $\mathcal{N}(A_G)$ represents the solution set to the system of equations
    \begin{align}
        \label{aux_6} u(z_j) + u(y_{j+a}) + u(y_{j+b}) &= 0 \qquad (j \in \Z_n),\\
        \label{aux_7} u(z_j) + u(x_{j-a}) + u(x_{j-b}) &= 0 \qquad (j \in \Z_n),\\
        \label{aux_5} u(x_j) + u(y_j) + u(z_{j+\frac{n}{2}}) &= 0 \qquad (j \in \Z_n) 
    \end{align}
in $u \in \mathbb{R}^{V_G}$. Suppose that a fixed vector $u \in \mathbb{R}^{V_G}$ is indeed a solution vector to the aforementioned system. From Equation \eqref{aux_5}  we quickly see that
    \begin{equation}\label{aux_8}
        u(z_j) = -u(x_{j+\frac{n}{2}}) - u(y_{j+\frac{n}{2}})
    \end{equation}
    holds for each $j \in \Z_n$. By plugging in Equation \eqref{aux_8}  into Equation \eqref{aux_6}, we further obtain
    \[
        -u(x_{j+\frac{n}{2}}) - u(y_{j+\frac{n}{2}}) + u(y_{j+a}) + u(y_{j+b}) = 0,
    \]
    which means that
    \begin{equation}\label{aux_9}
        u(x_j) = u(y_{j + \frac{n}{2} + a}) + u(y_{j+\frac{n}{2}+b}) - u(y_j)
    \end{equation}
is true for each $j \in \Z_n$. It is now possible to plug in Equations    \eqref{aux_8}  and 
\eqref{aux_9}  into Equation~\eqref{aux_7}  in order to get
    \begin{align*}
        \big(-u(y_{j+a}) - u(y_{j+b}) &+ u(y_{j+\frac{n}{2}}) - u(y_{j+\frac{n}{2}}) \big) \\
        &+ \left( u(y_{j + \frac{n}{2}}) + u(y_{j+\frac{n}{2} + b - a}) - u(y_{j-a}) \right) \\
        &+ \left( u(y_{j + \frac{n}{2} + a - b}) + u(y_{j + \frac{n}{2}}) - u(y_{j-b}) \right) = 0,
    \end{align*}
    i.e.,
    \begin{align}\label{aux_10}\begin{split}
        -u(y_{j+a}) - u(y_{j+b}) &- u(y_{j-a}) - u(y_{j-b})\\
        &+ u(y_{j+\frac{n}{2} + b - a}) + u(y_{j+\frac{n}{2} + a - b}) + 2 u(y_{j + \frac{n}{2}}) = 0,
    \end{split}\end{align}
    for each $j \in \Z_n$.
    
    Now, Equation \eqref{aux_10}  can be thought of as a system of equations in $u \in \mathbb{R}^Y$ and it is not difficult to establish that a vector $u$ represents its solution if and only if $\begin{bmatrix} u(y_0) & u(y_1) & \cdots & u(y_{n-1}) \end{bmatrix}^\intercal$ is a null space vector of the corresponding circulant matrix $C \in \mathbb{R}^{n \times n}$. By implementing Equation \eqref{circulant_eigenvalues}, we see that the eigenvalues of $C$ are obtained by the expression
    \begin{equation}\label{aux_11}
        -\zeta^a - \zeta^b - \zeta^{-a} - \zeta^{-b} + \zeta^{\frac{n}{2} + b - a} + \zeta^{\frac{n}{2} + a - b} + 2 \zeta^{\frac{n}{2}},
    \end{equation}
    as $\zeta$ ranges through the $n$-th roots of unity. We will now prove that the graph $G$ is a nut graph if and only if $C$ is of nullity one. 
    
    First of all, it is clear that the scenario $\eta(C) = 0$ is certainly not possible due to the fact that plugging in $\zeta = 1$ into Equation \eqref{aux_11}  yields the value zero. If $\eta(C) \ge 2$, then there exist two linearly independent solutions to Equation \eqref{aux_10}. However, if we apply Equations \eqref{aux_8} and \eqref{aux_9}  to these solutions, we further obtain two linearly independent solution vectors to the starting system determined by Equations \eqref{aux_6}, \eqref{aux_7}  and  \eqref{aux_5}. Thus, we have $\eta(G) \ge 2$, which means that $G$ is not a nut graph. On the other hand, if $\eta(C) = 1$, it becomes sufficient to notice that the solution set of Equation \eqref{aux_10}  contains the vectors $u \in \mathbb{R}^Y$ such that $u(y_j)$ is constant for each $j \in \Z_n$. The condition $\eta(C) = 1$ guarantees that these vectors are actually the only solutions to Equation~\eqref{aux_10}. By implementing Equations \eqref{aux_9}  and  \eqref{aux_8}, it can be swiftly seen that if $u \in \mathcal{N}(A_G)$, then this vector must be of the form
    \begin{align*}
        u(x_j) &= \beta, \quad & u(y_j) &= \beta, \quad & u(z_j) &= -2\beta \qquad && (j \in \Z_n),
    \end{align*}
    for some $\beta \in \mathbb{R}$. The converse is also trivial to check, which immediately tells us that $\eta(G) = 1$ and that $\mathcal{N}(A_G)$ contains a full vector. Hence, $G$ is a nut graph.

    We have proved that $G$ is a nut graph if and only if Equation \eqref{aux_11}  becomes zero for only a single value of $\zeta \in \mathbb{C}$ among the $n$-th roots of unity. Since $\zeta = 1$ necessarily yields the value zero, it is clear that this condition is equivalent to Equation \eqref{aux_11} being nonzero for each $\zeta \in \mathbb{C}, \, \zeta^n = 1, \, \zeta \neq 1$. For each such $\zeta$, it is sufficient to notice that
$$
-\zeta^a - \zeta^b - \zeta^{-a} - \zeta^{-b} + \zeta^{\frac{n}{2} + b - a} + \zeta^{\frac{n}{2} + a - b} + 2 \zeta^{\frac{n}{2}} = 0$$
if and only if
$$
 -\zeta^{a + b} \left( -\zeta^a - \zeta^b - \zeta^{-a} - \zeta^{-b} + \zeta^{\frac{n}{2} + b - a} + \zeta^{\frac{n}{2} + a - b} + 2 \zeta^{\frac{n}{2}} \right) = 0
$$
if and only if
$$
\zeta^{2a + b} + \zeta^{a + 2b} + \zeta^a + \zeta^b - \zeta^{\frac{n}{2} + 2a} - \zeta^{\frac{n}{2} + 2b} - 2\zeta^{\frac{n}{2} + a + b}=0,$$
%
in order to complete the proof of the lemma.
\end{proof}

We are now in position to apply the auxiliary Lemma \ref{type_b_polynomial_lemma} in order to give a relatively short proof of Theorem \ref{type_b_th}.

\bigskip\noindent
\emph{Proof of Theorem \ref{type_b_th}}.\quad
Let $G$ be a graph representable as $T_1(n, a, b)$. By virtue of Lemma \ref{type_b_polynomial_lemma}, we see that $G$ is a nut graph if and only if the polynomial given in Equation (\ref{aux_12}) contains only the root $1$ among all the $n$-th roots of unity. However, for any $x \in \mathbb{C}, \, x^n = 1$, it is possible to notice that
\begin{align*}
x^{2a + b} + x^{a + 2b} + x^a + x^b - x^{\frac{n}{2} + 2a} &- x^{\frac{n}{2} + 2b} - 2x^{\frac{n}{2} + a + b} =\\
&= x^a \, (x^{b-a} + 1)(x^{\frac{n}{2} + a} - 1)(x^{\frac{n}{2} + b} - 1) .
\end{align*}
From here, it is easy to see that $G$ is a nut graph if and only if the equation $x^{b-a} = -1$ has no solutions in $x \in \mathbb{C}, \, x^n = 1$, while the equations $x^{\frac{n}{2} + a} = 1$ and $x^{\frac{n}{2} + b} = 1$ each have a single solution in $x \in \mathbb{C}, \, x^n = 1$, namely the value $1$.

If we put $x = e^{\frac{2t \pi}{n} i}$ for a uniquely defined $t \in \mathbb{N}_0, \, 0 \le t < n$, the equation $x^{b-a} = -1$ in $x \in \mathbb{C}, \, x^n = 1$ becomes equivalent to the equation
\begin{equation}\label{aux_13}
    t(b-a) \equiv_n \frac{n}{2}
\end{equation}
in $t \in \mathbb{N}_0, \, 0 \le t < n$. Thus, the two equations will have the same number of solutions. However, Equation (\ref{aux_13}) represents a linear congruence equation, which means that it contains a solution if and only if $\gcd(n, b-a) \mid \frac{n}{2}$ (see, for example, \cite[p.~170, Theorem 5.14]{Tattersall}). Furthermore, it is simple to see that $\gcd(n, b-a) \mid \frac{n}{2}$ is equivalent to $v_2(b - a) < v_2(n)$. Thus, the starting equation has no solutions if and only if $v_2(b - a) \ge v_2(n)$.

In an analogous fashion, we may analyze the remaining two equations and discover that they both necessarily contain a solution. Moreover, the equation $x^{\frac{n}{2} + a} = 1$ has $\gcd(n, \frac{n}{2} + a)$ distinct solutions in total, while the number of solutions of $x^{\frac{n}{2} + b} = 1$ is equal to $\gcd(n, \frac{n}{2} + b)$. Taking everything into consideration, we conclude that $G$ is a nut graph if and only if the following conditions hold:
\begin{itemize}
    \item $\gcd(n, \frac{n}{2} + a) = 1$;
    \item $\gcd(n, \frac{n}{2} + b) = 1$;
    \item $v_2(b - a) \ge v_2(n)$.
\end{itemize}

We will now show that $\gcd(n, \frac{n}{2} + a) = 1$ is equivalent to the conjunction of $\gcd(\frac{n}{2}, a) = 1$ and $ a \not\equiv_2 \frac{n}{2}$. First of all, it is not difficult to realize that $\gcd(\frac{n}{2}, a) \mid \gcd(n, \frac{n}{2} + a)$, hence if $\gcd(\frac{n}{2}, a) > 1$, then $\gcd(n, \frac{n}{2} + a) > 1$ as well. Besides that, if $a \equiv_2 \frac{n}{2}$, then $n$ and $\frac{n}{2} + a$ are both surely even, which means that $\gcd(n, \frac{n}{2} + a) \ge 2$. These observations show that
$\gcd\left(n, \frac{n}{2} + a\right) = 1$ implies that $\gcd\left(\frac{n}{2}, a\right) = 1$ and $a \not\equiv_2 \frac{n}{2}$.

Now, suppose that $\gcd(\frac{n}{2}, a) = 1$ and $a \not\equiv_2 \frac{n}{2}$ do both hold. If we put $\beta = \gcd(n, \frac{n}{2} + a)$, it is clear that $\beta \mid n$ and $\beta \mid n + 2a$, which directly implies that $\beta \mid 2a$ as well, hence $\beta \mid \gcd(n, 2a)$. Thus, we obtain $\beta \mid 2\gcd(\frac{n}{2}, a)$, which is only possible if $\beta \in \{1, 2 \}$. Bearing in mind that $a \not\equiv_2 \frac{n}{2}$, it is clear that $\frac{n}{2} + a$ is odd, which implies $\beta \neq 2$, hence $\gcd(n, \frac{n}{2} + a) = 1$, as desired.

Thus, we have demonstrated that $\gcd(n, \frac{n}{2} + a) = 1$ is true if and only if $\gcd(\frac{n}{2}, a) = 1$ and $a \not\equiv_2 \frac{n}{2}$ are, and it is analogous to prove that $\gcd(n, \frac{n}{2} + b) = 1$ is equivalent to the conjunction of $\gcd(\frac{n}{2}, b) = 1$ and $b \not\equiv_2 \frac{n}{2}$. From here, we quickly conclude that $G$ being a nut graph is indeed equivalent to the conjunction of conditions \textbf{(i)}, \textbf{(ii)} and \textbf{(iii)}. \hfill\qed


\begin{remark}
    For any $t \in \mathbb{N}$ such that $\gcd(n, t) = 1$, it is not difficult to show that $T_1(n, a, b), \, 0 \le a, b < n, \, a \neq b$ is surely isomorphic to $T_1(n, ta \bmod n, tb \bmod n)$. Bearing this in mind, we may impose $n, a, b$ parameter conditions precisely required for $T_1(n, a, b)$ to be a nut graph which are stricter than those given in Theorem~\ref{type_b_th}. If $4 \mid n$, then a tricirculant graph of type 1 is a nut graph if and only if it is representable as $T_1(n, 1, b)$ where $3 \le b < n$, $\gcd(n, b) = 1$ and $v_2(b - 1) \ge v_2(n)$. Similarly, if $4 \nmid n$, then a tricirculant graph of type 1 is a nut graph if and only if it is representable as $T_1(n, 2, b)$ where $4 \le b < n$ and $\gcd(n, b) = 2$.
\end{remark}

\section{Tricirculant graphs of type 4}\label{sc_type_a}

In this final section, we shall determine all the nut graphs among the cubic tricirculant graphs of type $4$, thereby completing the proof of Theorem \ref{main_theorem}. The aforementioned result is disclosed within the following theorem.

\begin{theorem}\label{type_a_th}
    An arbitrary graph representable as $T_4(n, a, b)$ where $n$ is even and $1 \le a, b < \frac{n}{2}$, is a nut graph if and only if the following conditions hold:
    \begin{enumerate}[label=\textbf{(\roman*)}]
        \item $\gcd\left( \frac{n}{2}, a, b \right) = 1$;
        \item if $4 \nmid n$, then at least one of $a, b$ is even;
        \item if $4 \mid n$, then $a$ and $b$ are of different parities;
        \item if $10 \mid n$, then at least one of $a, b, a-b, a+b$ is divisible by five.
    \end{enumerate}
\end{theorem}

We will now demonstrate how the problem of testing whether a $T_4(n, a, b)$ graph is a nut graph can be transformed to a number theory problem, in a similar manner as it was done in Lemma \ref{type_b_polynomial_lemma}. The corresponding result is given in the next lemma.

\begin{lemma}\label{type_a_polynomial_lemma}
    A graph representable as $T_4(n, a, b)$ is a nut graph if and only if the $\mathbb{Z}[x]$ polynomial
    \begin{equation}\label{aux_14}
         x^{2a + b} + x^{a + 2b} + x^a + x^b - x^{\frac{n}{2} + 2a + 2b} - x^{\frac{n}{2} + 2a} - x^{\frac{n}{2} + 2b} - x^{\frac{n}{2}}
    \end{equation}
    has no $n$-th roots of unity among its roots, besides $1$.
\end{lemma}
\begin{proof}
    Let $G$ be a given graph that is representable as $T_4(n, a, b)$. The very definition of the type 4 graphs dictates that $\mathcal{N}(A_G)$ represents the solution set to the system of equations
    \begin{align}
        \label{aux_15} u(x_j) + u(y_j) + u(z_{j+\frac{n}{2}}) &= 0 \qquad (j \in \Z_n),\\
        \label{aux_16} u(z_j) + u(y_{j+b}) + u(y_{j-b}) &= 0 \qquad (j \in \Z_n),\\
        \label{aux_17} u(z_j) + u(x_{j+a}) + u(x_{j-a}) &= 0 \qquad (j \in \Z_n)
    \end{align}
    in $u \in \mathbb{R}^{V_G}$. If we suppose that some vector $u \in \mathbb{R}^{V_G}$ is a solution vector to the given system, Equation (\ref{aux_15}) lets us immediately obtain that
    \begin{equation}\label{aux_18}
        u(z_j) = -u(x_{j+\frac{n}{2}}) - u(y_{j+\frac{n}{2}})
    \end{equation}
    is true for each $j \in \Z_n$. Furthermore, we can plug in Equation (\ref{aux_18}) into Equation~(\ref{aux_16}) in order to get
    \[
        -u(x_{j+\frac{n}{2}}) - u(y_{j+\frac{n}{2}}) + u(y_{j+b}) + u(y_{j-b}) = 0,
    \]
    i.e.,
    \begin{equation}\label{aux_19}
        u(x_j) = u(y_{j+\frac{n}{2}+b}) + u(y_{j+\frac{n}{2}-b}) - u(y_j)
    \end{equation}
    for each $j \in \Z_n$. Now, by plugging in Equations (\ref{aux_18}) and (\ref{aux_19}) into Equation (\ref{aux_17}), we may conclude that
    \begin{align*}
        \big(-u(y_{j+b}) - u(y_{j-b}) &+ u(y_{j+\frac{n}{2}}) - u(y_{j+\frac{n}{2}}) \big) \\
        &+ \left( u(y_{j + \frac{n}{2} + b + a}) + u(y_{j + \frac{n}{2} - b + a}) - u(y_{j+a}) \right) \\
        &+ \left( u(y_{j + \frac{n}{2} + b - a}) + u(y_{j + \frac{n}{2} - b - a}) - u(y_{j-a}) \right) = 0,
    \end{align*}
    which immediately gives
    \begin{align}\label{aux_20}\begin{split}
        -u(y_{j+a}) - u(y_{j+b}) &- u(y_{j-a}) - u(y_{j-b}) + u(y_{j+\frac{n}{2} + a + b})\\
        &+ u(y_{j+\frac{n}{2} + a - b}) + u(y_{j + \frac{n}{2} - a + b}) + u(y_{j + \frac{n}{2} - a - b}) = 0,
    \end{split}\end{align}
    for each $j \in \Z_n$.

    By thinking of Equation (\ref{aux_20}) as a system of equations in $u \in \mathbb{R}^Y$, we are now able to implement the same circulant matrix null space strategy from Lemma \ref{type_b_polynomial_lemma}. In other words, by using an analogous proof method which we choose to leave out for the sake of brevity, we may conclude that the graph $G$ is a nut graph if and only if the expression
    \begin{equation}\label{aux_21}
        -\zeta^a - \zeta^b - \zeta^{-a} - \zeta^{-b} + \zeta^{\frac{n}{2} + a + b} + \zeta^{\frac{n}{2} + a - b} + \zeta^{\frac{n}{2} - a + b} + \zeta^{\frac{n}{2} - a - b}
    \end{equation}
    yields zero only for a single value of $\zeta \in \mathbb{C}, \, \zeta^n = 1$. Given the fact that Equation~(\ref{aux_21}) surely gives zero when we plug in $\zeta = 1$, the said condition must be equivalent to Equation~(\ref{aux_21}) not giving zero for any other $n$-th root of unity. However, for each $\zeta \in \mathbb{C}, \, \zeta^n = 1$, we have
$$-\zeta^a - \zeta^b - \zeta^{-a} - \zeta^{-b} + \zeta^{\frac{n}{2} + a + b} + \zeta^{\frac{n}{2} + a - b} + \zeta^{\frac{n}{2} - a + b} + \zeta^{\frac{n}{2} - a - b}  = 0$$
if and only if
$$ -\zeta^{a + b} \left( -\zeta^a - \zeta^b - \zeta^{-a} - \zeta^{-b} + \zeta^{\frac{n}{2} + a + b} + \zeta^{\frac{n}{2} + a - b} + \zeta^{\frac{n}{2} - a + b} + \zeta^{\frac{n}{2} - a - b} \right)  = 0$$
if and only if
$$\zeta^{2a + b} + \zeta^{a + 2b} + \zeta^a + \zeta^b - \zeta^{\frac{n}{2} + 2a + 2b} - \zeta^{\frac{n}{2} + 2a} - \zeta^{\frac{n}{2} + 2b} - \zeta^{\frac{n}{2}}  = 0,$$
which promptly leads us to the desired lemma statement.
\end{proof}

Unfortunately, the Equation (\ref{aux_14}) polynomial does not adhere to a factorization analogous to the one applied on the Equation (\ref{aux_12}) polynomial during the proof of Theorem~\ref{type_b_th}. For this reason, we shall use an entirely different strategy in order to complete the proof of Theorem \ref{type_a_th}. To begin, we define the following two auxiliary polynomials
\begin{align*}
    Q_{a, b}(x) &= x^{2a + b} + x^{a + 2b} + x^a + x^b - x^{2a + 2b} - x^{2a} - x^{2b} - 1,\\
    R_{a, b}(x) &= x^{2a + b} + x^{a + 2b} + x^a + x^b + x^{2a + 2b} + x^{2a} + x^{2b} + 1,
\end{align*}
for each $a, b \in \mathbb{N}$. It now becomes convenient to disclose the following brief reformulation of the polynomial problem obtained as a result of Lemma \ref{type_a_polynomial_lemma}.

\begin{lemma}\label{type_a_cyclotomic_lemma}
    A graph representable as $T_4(n, a, b)$ is a nut graph if and only if the following two conditions hold:
    \begin{enumerate}[label=\textbf{(\roman*)}]
        \item $\Phi_f(x) \nmid Q_{a, b}(x)$ for every $f \in \mathbb{N}, \, f \ge 2$ such that $f \mid \frac{n}{2}$;
        \item $\Phi_f(x) \nmid R_{a, b}(x)$ for every $f \in \mathbb{N}$ such that $f \mid n$ and $2 \nmid \frac{n}{f}$.
    \end{enumerate}
\end{lemma}
\begin{proof}
    Let $G$ be a given graph that is representable as $T_4(n, a, b)$. By virtue of Lemma \ref{type_a_polynomial_lemma}, it can be immediately seen that $G$ is a nut graph if and only if the following two conditions hold:
    \begin{itemize}
        \item $Q_{a, b}(x)$ does not contain a root $\zeta$ among all the $n$-th roots of unity such that $\zeta^\frac{n}{2} = 1$, besides $1$;
        \item $R_{a, b}(x)$ does not contain a root $\zeta$ among all the $n$-th roots of unity such that $\zeta^\frac{n}{2} = -1$.
    \end{itemize}
    A complex number $\zeta \in \mathbb{C}$ is an $n$-th root of unity and satisfies $\zeta^\frac{n}{2} = 1$ if and only if it is an $f$-th primitive root of unity for some $f\in\mathbb{N}$ such that $f \mid \frac{n}{2}$. From here, it is not difficult to realize that $Q_{a, b}(x)$ does not contain a root among all such $n$-th roots of unity, besides $1$, if and only if it is not divisible by any cyclotomic polynomial $\Phi_f(x)$ where $f \ge 2$ and $f \mid \frac{n}{2}$. In a similar fashion, it can be noticed that some $\zeta \in \mathbb{C}$ is an $n$-th root of unity and satisfies $\zeta^\frac{n}{2} = -1$ if and only if this number is an $f$-th primitive root of unity for some $f \in \mathbb{N}$ such that $f \mid n$ and $2 \nmid \frac{n}{f}$. However, this means that $R_{a, b}(x)$ does not contain such a root if and only if it is not divisible by any polynomial $\Phi_f(x)$ where $f \mid n$ and $2 \nmid \frac{n}{f}$.
\end{proof}

It is worth pointing out that, for any $f \in \mathbb{N}$, $\Phi_f(x)$ divides a given polynomial $V(x) \in \mathbb{Q}[x]$ if and only if it divides any other polynomial which can be obtained from $V(x)$ by adding or subtracting any multiple of $f$ from the powers of its terms. Bearing this in mind, our next step in proving Theorem \ref{type_a_th} will be to demonstrate the validity of the easier underlying implication of the required equivalence. In other words, we will show that each graph not satisfying the stated conditions is certainly not a nut graph. This result is given in the next lemma.

\begin{lemma}\label{type_a_easy_direction}
    A graph representable as $T_4(n, a, b)$ which does not satisfy the conditions stated in Theorem \ref{type_a_th} is surely not a nut graph.
\end{lemma}
\begin{proof}
    Let $G$ be a given graph which is representable as $T_4(n, a, b)$ and suppose that this graph does not satisfy all four conditions disclosed in Theorem \ref{type_a_th}. In order to make the proof more concise, we will split it into four cases depending on which of the four conditions is not satisfied. In all the cases, we shall demonstrate that $G$ cannot be a nut graph.

    \bigskip\noindent
    \textbf{Case  1:} $\gcd(\frac{n}{2}, a, b) = 1$ does not hold.\quad
    Let $\beta = \gcd(\frac{n}{2}, a, b) \ge 2$. Given the fact that both $a$ and $b$ are divisible by $\beta$, it is straightforward to notice that $Q_{a, b}(x)$ must contain each $\beta$-th root of unity among its roots, hence $\Phi_\beta(x) \mid Q_{a, b}(x)$. However, since $\beta \ge 2$ and $\beta \mid \frac{n}{2}$, condition \textbf{(i)} from Lemma \ref{type_a_cyclotomic_lemma} dictates that $G$ cannot be a nut graph, as desired.

    \bigskip\noindent
    \textbf{Case  2:} $4 \nmid n \implies 2 \mid ab$ does not hold.\quad
    Suppose that $4 \nmid n$, while $2 \nmid a$ and $2 \nmid b$. It is obvious that $2 \mid n$ and $\frac{n}{2}$ is odd. Also, it can be swiftly noticed that $-1$ must be a root of $R_{a, b}(x)$, which means that $\Phi_2(x) \mid R_{a, b}(x)$. Thus, condition \textbf{(ii)} from Lemma \ref{type_a_cyclotomic_lemma} fails, which implies that the graph $G$ cannot be a nut graph.
    
    \bigskip\noindent
    \textbf{Case 3:} $4 \mid n \implies a \not\equiv_2 b$ does not hold.\quad Now, suppose that $4 \mid n$, while $a$ and $b$ are of the same parity. If $a$ and $b$ are both even, then it is easy to see that $-1$ is a root of $Q_{a, b}(x)$. Thus, $\Phi_2(x) \mid Q_{a, b}(x)$ must hold, while $2 \mid \frac{n}{2}$. By implementing condition \textbf{(i)} from Lemma \ref{type_a_cyclotomic_lemma}, we get that $G$ is not a nut graph, as desired.

    If $a$ and $b$ are both odd, then it becomes convenient to split the problem into two further subcases depending on whether $8 \mid n$ or $n \equiv_8 4$. If $n \equiv_8 4$, then $4 \mid n$ and $2 \nmid \frac{n}{4}$, together with
    \[
        R_{a, b}(i)  = i^{2a + b} + i^{a+2b} + i^a + i^b + i^{2a + 2b} + i^{2a} + i^{2b} + 1 = (i^{a + b} + 1)(i^a + i^b) .
    \]
    If $a \equiv_4 b$, then $i^{a + b} + 1 = 0$, and if $a \not\equiv_4 b$, then $i^a + i^b = 0$. Either way, we conclude that $R_{a, b}(i) = 0$, which further implies $\Phi_4(x) \mid R_{a, b}(x)$. According to condition \textbf{(ii)} from Lemma \ref{type_a_cyclotomic_lemma}, the graph $G$ is not a nut graph. Finally, if $8 \mid n$, we may observe that $Q_{a, b}(i) = 0$ by using the same computational strategy as done in the previous subcase. This directly gives $\Phi_4(x) \mid Q_{a, b}(x)$. Bearing in mind that $4 \mid \frac{n}{2}$, condition~\textbf{(i)} from Lemma~\ref{type_a_cyclotomic_lemma} tells us that $G$ is not a nut graph.

    \bigskip\noindent
    \textbf{Case  4:} $10 \mid n \implies 5 \mid ab(a-b)(a+b)$  does not hold.\quad
    Finally, suppose that $10 \mid n$ and that $5 \nmid ab(a-b)(a+b)$. It is not difficult to establish that the latter condition is equivalent to the disjuction of $a \bmod 5 \in \{ 1, 4\} \land b \bmod 5 \in \{ 2, 3 \}$ and \linebreak $b \bmod 5 \in \{ 1, 4\} \land a \bmod 5 \in \{ 2, 3 \}$. Given the fact that $5 \mid \frac{n}{2}$, Lemma \ref{type_a_cyclotomic_lemma} claims that in order to show that $G$ is not a nut graph, it is sufficient to prove that $\Phi_5(x) \mid Q_{a, b}(x)$, and this is exactly what we shall do. Since $Q_{a, b}(x) = Q_{b, a}(x)$, we may assume that $a \bmod 5 \in \{ 1, 4\}$ and $b \bmod 5 \in \{ 2, 3 \}$, without loss of generality. From here, it is straightforward to notice that either $b \equiv_5 2a$ or $b \equiv_5 3a$.

    If we have $b \equiv_5 2a$, then by plugging in the said modular equality we swiftly notice that
    \begin{alignat*}{2}
        && \Phi_5(x) &\mid Q_{a, b}(x)\\
        \iff \quad && \Phi_5(x) &\mid x^{4a} + x^{5a} + x^a + x^{2a} - x^{6a} - x^{2a} - x^{4a} - 1\\
        \iff \quad && \Phi_5(x) &\mid x^{4a} + 1 + x^a + x^{2a} - x^{a} - x^{2a} - x^{4a} - 1\\
        \iff \quad && \Phi_5(x) &\mid 0,
    \end{alignat*}
    hence $\Phi_5(x) \mid Q_{a, b}(x)$ is indeed true. Similarly, if $b \equiv_5 3a$, then by using the same strategy we obtain
    \begin{alignat*}{2}
        && \Phi_5(x) &\mid Q_{a, b}(x)\\
        \iff \quad && \Phi_5(x) &\mid x^{5a} + x^{7a} + x^a + x^{3a} - x^{8a} - x^{2a} - x^{6a} - 1\\
        \iff \quad && \Phi_5(x) &\mid 1 + x^{2a} + x^a + x^{3a} - x^{3a} - x^{2a} - x^{a} - 1\\
        \iff \quad && \Phi_5(x) &\mid 0,
    \end{alignat*}
    which directly implies $\Phi_5(x) \mid Q_{a, b}(x)$ once again.
\end{proof}

Before we proceed with the proof of Theorem \ref{type_a_th}, we will need another short auxiliary folklore lemma.

\begin{lemma}\label{modulo_div_lemma}
    Let $V(x), W(x) \in \mathbb{Q}[x], \, W(x) \not\equiv 0$ be two polynomials such that $W(x) \mid V(x)$ and the powers of all the nonzero terms of $W(x)$ are divisible by some $\beta \in \mathbb{N}$. 
For any $j \in\{0,\dots, \beta-1\}$, if we use $V^{(\beta, j)}(x)$ to denote the polynomial composed of all the terms of $V(x)$ whose powers are congruent to $j$ modulo $\beta$, we then have $W(x) \mid V^{(\beta, j)}(x)$.
\end{lemma}
\begin{proof}
    Since $W(x) \mid V(x)$, we may write $V(x) = W(x) \, U(x)$ for some polynomial $U(x) \in \mathbb{Q}[x]$. For each $j \in\{0,\dots, \beta-1\}$, let $U^{(\beta, j)}(x)$ denote the polynomial comprising all the terms of $U(x)$ whose powers are congruent to $j$ modulo $\beta$. It now becomes simple to notice that $V^{(\beta, j)}(x) = W(x) \, U^{(\beta, j)}(x)$ must be true for each $j \in\{0,\dots, \beta-1\}$. The lemma statement follows directly from here.
\end{proof}

We shall now extensively implement Lemma \ref{modulo_div_lemma} in order to obtain a series of auxiliary lemmas regarding the divisibility of $Q_{a, b}(x)$ and $R_{a, b}(x)$ by polynomials whose nonzero terms have powers divisible by a common prime.

\begin{lemma}\label{complicated_lemma_p}
    If a given $W(x) \in \mathbb{Q}[x]$ has at least two nonzero terms and all of its nonzero terms have powers divisible by a prime number $p \ge 7$, then $W(x) \mid Q_{a, b}(x)$ and $W(x) \mid R_{a, b}(x)$ both imply $p \mid a, b$.
\end{lemma}
\begin{proof}
    First of all, if all the numbers $2a + b, a + 2b, a, b, 2a + 2b, 2a, 2b$ are not divisible by $p$, we then obtain
    \begin{align*}
        Q_{a, b}^{(p, 0)}(x) = -1, \quad R_{a, b}^{(p, 0)}(x) = 1 .
    \end{align*}
    From here we have that $W(x) \nmid Q_{a, b}^{(p, 0)}(x), R_{a, b}^{(p, 0)}(x)$, hence $W(x) \nmid Q_{a, b}(x), R_{a, b}(x)$, by virtue of Lemma \ref{modulo_div_lemma}. Thus, there is nothing left to discuss in this scenario. We will now suppose that at least one number from $2a + b, a + 2b, a, b, 2a + 2b, 2a, 2b$ is divisible by $p$. In order to make the proof easier to follow, we shall split the problem into seven corresponding cases.

    \bigskip\noindent
    \textbf{Case 1:} $p \mid 2a + b$.\quad
    Here, we have that $b \equiv_p -2a$, which implies
    \begin{align*}
        2a + b &\equiv_p 0, & a + 2b &\equiv_p -3a, & a &\equiv_p a, & b &\equiv_p -2a,\\
        2a + 2b &\equiv_p -2a, & 2a &\equiv_p 2a, & 2b &\equiv_p -4a, & 0 &\equiv_p 0.
    \end{align*}
    It is now easy to see that we obtain two further possibilities:
    \begin{itemize}
        \item $a$ has a unique remainder modulo $p$ within the set $\{2a + b, a + 2b, a, b, \linebreak 2a + 2b, 2a, 2b, 0 \}$;
        \item at least one number from the set $\{ a, 3a, 4a, 5a \}$ is divisible by $p$.
    \end{itemize}
    In the former scenario, Lemma \ref{modulo_div_lemma} dictates that $W(x) \mid Q_{a, b}(x)$ or $W(x) \mid R_{a, b}(x)$ would both imply $W(x) \mid x^a$, which is not possible. Thus, $W(x) \nmid Q_{a, b}(x), R_{a, b}(x)$ and the lemma statement holds. In the latter scenario, it is not difficult to conclude that each subcase leads to $p \mid a$. From here, it immediately follows that $p \mid a, b$, as desired. 

    \bigskip\noindent
    \textbf{Case 2:} $p \mid a + 2b$.\quad
    This case can be proved in an entirely analogous manner as case 1.
    
    \bigskip\noindent
    \textbf{Case 3:} $p \mid a$.\quad
    In this case, we get $a \equiv_p 0$, which means that
    \begin{align*}
        2a + b &\equiv_p b, & a + 2b &\equiv_p 2b, & a &\equiv_p 0, & b &\equiv_p b,\\
        2a + 2b &\equiv_p 2b, & 2a &\equiv_p 0, & 2b &\equiv_p 2b, & 0 &\equiv_p 0.
    \end{align*}
    We now get two further possibilities:
    \begin{itemize}
        \item $p \mid b$ or $p \mid 2b$;
        \item the numbers $\{2a + b, a + 2b, a, b, 2a + 2b, 2a, 2b, 0 \}$ can be partitioned as $\{ \{2a + b, b \}, \{a + 2b, 2a + 2b, 2b \}, \{ a, 2a, 0 \} \}$ according to their remainder modulo $p$.
    \end{itemize}
    In the former scenario, we certainly have that $p \mid a, b$, hence there is nothing more to discuss. In the later scenario, Lemma \ref{modulo_div_lemma} allows us to swiftly conclude that $W(x) \mid Q_{a, b}(x)$ implies
    \begin{align}
        \label{aux_22} W(x) &\mid x^{2a + b} + x^b,\\
        \label{aux_23} W(x) &\mid x^{a + 2b} - x^{2a + 2b} - x^{2b},
    \end{align}
    which is not possible due to the fact that Equations (\ref{aux_22}) and (\ref{aux_23}) together give
    \begin{alignat*}{2}
        && W(x) &\mid (x^{a + 2b} - x^{2a + 2b} - x^{2b}) + x^b \, (x^{2a + b} + x^b)\\
        \implies \quad && W(x) &\mid x^{a + 2b} .
    \end{alignat*}
    In a similar fashion, $W(x) \mid R_{a, b}(x)$ would lead to
    \begin{align}
        \label{aux_24} W(x) &\mid x^{2a + b} + x^b,\\
        \label{aux_25} W(x) &\mid x^{a + 2b} + x^{2a + 2b} + x^{2b},
    \end{align}
    However, this is again impossible since Equations (\ref{aux_24}) and (\ref{aux_25}) imply
    \begin{alignat*}{2}
        && W(x) &\mid (x^{a + 2b} + x^{2a + 2b} + x^{2b}) - x^b \, (x^{2a + b} + x^b)\\
        \implies \quad && W(x) &\mid x^{a + 2b} .
    \end{alignat*}

    \bigskip\noindent
    \textbf{Case 4:} $p \mid b$.\quad
    This case can be proved in an entirely analogous manner as case~3.

    \bigskip\noindent
    \textbf{Case 5:} $p \mid 2a + 2b$.\quad
    The condition $p \mid 2a + 2b$ gets down to $b \equiv_p -a$, which immediately leads to
    \begin{align*}
        2a + b &\equiv_p a, & a + 2b &\equiv_p -a, & a &\equiv_p a, & b &\equiv_p -a,\\
        2a + 2b &\equiv_p 0, & 2a &\equiv_p 2a, & 2b &\equiv_p -2a, & 0 &\equiv_p 0.
    \end{align*}
    From here, we reach two further possibilities:
    \begin{itemize}
        \item $2a$ has a unique remainder modulo $p$ within the set $\{2a + b, a + 2b, a, b, \linebreak 2a + 2b, 2a, 2b, 0 \}$;
        \item at least one number from the set $\{ a, 2a, 3a, 4a \}$ is divisible by $p$;
    \end{itemize}
    In the former scenario, Lemma \ref{modulo_div_lemma} states that $W(x) \mid Q_{a, b}(x)$ or $W(x) \mid R_{a, b}(x)$ would both imply $W(x) \mid x^{2a}$, which is clearly impossible. For this reason, we obtain $W(x) \nmid Q_{a, b}(x), R_{a, b}(x)$ as desired. In the latter scenario, it can be quickly noticed that $p \mid a$ must be true. This leads to $p \mid a, b$, which completes the proof.

    \bigskip\noindent
    \textbf{Case 6:} $p \mid 2a$.\quad
    The condition $p \mid 2a$ directly gives us $a \equiv_p 0$, which means that this case actually coincides with case 3.

    \bigskip\noindent
    \textbf{Case 7:} $p \mid 2b$.\quad
    This case can be proved in an entirely analogous manner as case~6.
\end{proof}

\begin{lemma}\label{complicated_lemma_5}
    If $W(x) \in \mathbb{Q}[x]$ is a polynomial with at least two nonzero terms such that all of its nonzero terms have powers divisible by five, then $W(x) \mid Q_{a, b}(x)$ and $W(x) \mid R_{a, b}(x)$ both imply $5 \mid a, b$ or $5 \nmid a, b, a+b, a-b$.
\end{lemma}
\begin{proof}
    To begin, if $5 \mid a, b$, then the lemma statement certainly holds and there is nothing more to discuss. Due to the fact that $Q_{a, b}(x) = Q_{b, a}(x)$, as well as $R_{a, b}(x) = R_{b, a}(x)$, we can assume, without loss of generality, that $5 \nmid b$. Now, it is not difficult to notice that $a$ must satisfy precisely one of the following five modular equalies: $a \equiv_5 0$, $a \equiv_5 b$, $a \equiv_5 2b$, $a \equiv_5 3b$, $a \equiv_5 4b$. We shall use this observation to complete the proof by splitting the problem into five corresponding cases.

    \begin{table}[h!t]
    {\footnotesize
    \begin{center}
    \begin{tabular}{rccccc}
    \toprule & $a \equiv_5 0$ & $a \equiv_5 b$ & $a \equiv_5 2b$ & $a \equiv_5 3b$ & $a \equiv_5 4b$\\
    \midrule
    $2a+b \equiv_5$ & $b$ & $3b$ & $0$ & $2b$ & $4b$\\
    $a+2b \equiv_5$ & $2b$ & $3b$ & $4b$ & $0$ & $b$\\
    $a \equiv_5$ & $0$ & $b$ & $2b$ & $3b$ & $4b$\\
    $b \equiv_5$ & $b$ & $b$ & $b$ & $b$ & $b$\\
    $2a+2b \equiv_5$ & $2b$ & $4b$ & $b$ & $3b$ & $0$\\
    $2a \equiv_5$ & $0$ & $2b$ & $4b$ & $b$ & $3b$\\
    $2b \equiv_5$ & $2b$ & $2b$ & $2b$ & $2b$ & $2b$\\
    $0 \equiv_5$ & $0$ & $0$ & $0$ & $0$ & $0$\\
    \bottomrule
    \end{tabular}
    \end{center}
    \caption{The powers of the $Q_{a, b}(x)$ and $R_{a, b}(x)$ nonzero terms modulo five.}
    \label{complicated_remainders_5}
    }
    \end{table}

    \bigskip\noindent
    \textbf{Case 1:} $a \equiv_5 0$.\quad
    In this scenario, Table \ref{complicated_remainders_5} and Lemma \ref{modulo_div_lemma} together dictate that $W(x) \mid Q_{a, b}(x)$ implies
    \begin{align*}
        W(x) &\mid x^{2a + b} + x^b,\\
        W(x) &\mid x^{a + 2b} - x^{2a + 2b} - x^{2b},\\
        W(x) &\mid x^a - x^{2a} - 1 ,
    \end{align*}
    while $W(x) \mid R_{a, b}(x)$ implies
    \begin{align*}
        W(x) &\mid x^{2a + b} + x^b,\\
        W(x) &\mid x^{a + 2b} + x^{2a + 2b} + x^{2b},\\
        W(x) &\mid x^a + x^{2a} + 1.
    \end{align*}
    From here onwards, the proof can be completed by using the same strategy as done so in case 3 from the proof of Lemma \ref{complicated_lemma_p}.

    \bigskip\noindent
    \textbf{Case 2:} $a \equiv_5 b$.\quad
    In this case, Table \ref{complicated_remainders_5} states that the element $0$ has a unique remainder modulo five within the set $\{2a + b, a + 2b, a, b, 2a + 2b, 2a, 2b, 0 \}$. For this reason, $W(x) \mid Q_{a, b}(x)$ and $W(x) \mid R_{a, b}(x)$ would both imply $W(x) \mid 1$, by virtue of Lemma \ref{modulo_div_lemma}. Thus, neither $W(x) \mid Q_{a, b}(x)$ nor $W(x) \mid R_{a, b}(x)$ can be true.

    \bigskip\noindent
    \textbf{Case 3:} $a \equiv_5 2b$.\quad
    Here, it is easy to check that $5 \nmid a, b, a + b, a - b$ and there is nothing more to discuss.

    \bigskip\noindent
    \textbf{Case 4:} $a \equiv_5 3b$.\quad
    In this situation, it is straightforward to notice that once again $5 \nmid a, b, a + b, a - b$, which means that the lemma statement does hold.

    \bigskip\noindent
    \textbf{Case 5:} $a \equiv_5 4b$.\quad
    Here, Table \ref{complicated_remainders_5} tells us that the element $2a$ has a unique remainder modulo five within the set $\{2a + b, a + 2b, a, b, 2a + 2b, 2a, 2b, 0 \}$. The rest of the proof can now be carried out analogously as it was done in case 2.
\end{proof}

\begin{lemma}\label{complicated_lemma_3}
    If $W(x) \in \mathbb{Q}[x]$ is a polynomial with at least two nonzero terms such that all of its nonzero terms have powers divisible by three, then $W(x) \mid Q_{a, b}(x)$ and $W(x) \mid R_{a, b}(x)$ both imply $3\mid a, b$.
\end{lemma}
\begin{proof}
    For starters, if $3 \mid a, b$, then the lemma statement holds and there is nothing left to prove. Given the fact that $Q_{a, b}(x) = Q_{b, a}(x)$ and $R_{a, b}(x) = R_{b, a}(x)$, we may assume, without loss of generality, that $3 \nmid b$. It is easy to establish that either $a \equiv_3 0$ or $a \equiv_3 b$ or $a \equiv_3 2b$ must hold. For this reason, it becomes convenient to split the problem into the three corresponding cases that arise from this observation.

    \begin{table}[h!t]
    {\footnotesize
    \begin{center}
    \begin{tabular}{rccc}
    \toprule & $a \equiv_3 0$ & $a \equiv_3 b$ & $a \equiv_3 2b$ \\
    \midrule
    $2a+b \equiv_3$ & $b$ & $0$ & $2b$\\
    $a+2b \equiv_3$ & $2b$ & $0$ & $b$\\
    $a \equiv_3$ & $0$ & $b$ & $2b$\\
    $b \equiv_3$ & $b$ & $b$ & $b$\\
    $2a+2b \equiv_3$ & $2b$ & $b$ & $0$\\
    $2a \equiv_3$ & $0$ & $2b$ & $b$\\
    $2b \equiv_3$ & $2b$ & $2b$ & $2b$\\
    $0 \equiv_3$ & $0$ & $0$ & $0$\\
    \bottomrule
    \end{tabular}
    \end{center}
    \caption{The powers of the $Q_{a, b}(x)$ and $R_{a, b}(x)$ nonzero terms modulo three.}
    \label{complicated_remainders_3}
    }
    \end{table}

    \bigskip\noindent
    \textbf{Case 1:} $a \equiv_3 0$.\quad
    In this case, Table \ref{complicated_remainders_3} and Lemma \ref{modulo_div_lemma} tell us that $W(x) \mid Q_{a, b}(x)$ implies
    \begin{align*}
        W(x) &\mid x^{2a + b} + x^b,\\
        W(x) &\mid x^{a + 2b} - x^{2a + 2b} - x^{2b},\\
        W(x) &\mid x^a - x^{2a} - 1 ,
    \end{align*}
    while $W(x) \mid R_{a, b}(x)$ implies
    \begin{align*}
        W(x) &\mid x^{2a + b} + x^b,\\
        W(x) &\mid x^{a + 2b} + x^{2a + 2b} + x^{2b},\\
        W(x) &\mid x^a + x^{2a} + 1.
    \end{align*}
    It is now evident that the given case can be resolved in the same manner as case~3 from Lemma \ref{complicated_lemma_p}.

    \bigskip\noindent
    \textbf{Case 2:} $a \equiv_3 b$.\quad
    If we suppose that $W(x) \mid Q_{a, b}(x)$, then Table \ref{complicated_remainders_3} and Lemma \ref{modulo_div_lemma} yield
    \begin{align*}
        W(x) &\mid x^{2a + b} + x^{a + 2b} - 1,\\
        W(x) &\mid x^a + x^b - x^{2a + 2b},\\
        W(x) &\mid x^{2a} + x^{2b},
    \end{align*}
    which promptly leads us to
    \begin{alignat*}{2}
        && W(x) &\mid (x^b - x^a)(x^a + x^b)^2 \, (x^{2a + b} + x^{a + 2b} - 1)\\
        && & \hspace{0.9cm} + (x^{2b} - x^{2a})(x^a + x^b - x^{2a + 2b})\\
        && & \hspace{0.9cm} + x^{a + b}(x^{2a} - x^{2b} + x^{a + b})(x^{2a} + x^{2b})\\
        \implies \quad && W(x) &\mid 2x^{2a + 4b} ,
    \end{alignat*}
    which is impossible. Thus, $W(x) \nmid Q_{a, b}(x)$. In a similar fashion, if we suppose that $W(x) \mid R_{a, b}(x)$, it follows that
    \begin{align*}
        W(x) &\mid x^{2a + b} + x^{a + 2b} + 1,\\
        W(x) &\mid x^a + x^b + x^{2a + 2b},\\
        W(x) &\mid x^{2a} + x^{2b} .
    \end{align*}
    Subsequently, we may obtain
    \begin{alignat*}{2}
        && W(x) &\mid (x^a - x^b)(x^a + x^b)^2 \, (x^{2a + b} + x^{a + 2b} + 1)\\
        && & \hspace{0.9cm} + (x^{2b} - x^{2a})(x^a + x^b + x^{2a + 2b})\\
        && & \hspace{0.9cm} + x^{a + b}(x^{2b} - x^{2a} + x^{a + b})(x^{2a} + x^{2b})\\
        \implies \quad && W(x) &\mid 2x^{4a + 2b} ,
    \end{alignat*}
    which is again not possible, hence $W(x) \nmid R_{a, b}(x)$.

    \bigskip\noindent
    \textbf{Case 3:} $a \equiv_3 2b$.\quad
    If we suppose that $W(x) \mid Q_{a, b}(x)$, then Table \ref{complicated_remainders_3} and Lemma~\ref{modulo_div_lemma} dictate that
    \begin{align*}
        W(x) &\mid x^{2a + b} + x^a - x^{2b},\\
        W(x) &\mid x^{a + 2b} + x^b - x^{2a},\\
        W(x) &\mid x^{2a + 2b} + 1 .
    \end{align*}
    However, in this scenario we may conclude that
    \begin{alignat*}{2}
        && W(x) &\mid (-1-x^{12b}+x^{a+7b}) (x^{2a + b} + x^a - x^{2b})\\
        && & \hspace{0.9cm} + (-x^b + x^{a+2b} + x^{a+11b} - x^{2a+6b})(x^{a + 2b} + x^b - x^{2a})\\
        && & \hspace{0.9cm} + (x^a + x^{a + 9b} - x^{2a + 4b})(x^{2a + 2b} + 1)\\
        \implies \quad && W(x) &\mid x^{14b} ,
    \end{alignat*}
    which is obviously impossible. Similarly, if we suppose that $W(x) \mid R_{a, b}(x)$, we immediately get
    \begin{align*}
        W(x) &\mid x^{2a + b} + x^a + x^{2b},\\
        W(x) &\mid x^{a + 2b} + x^b + x^{2a},\\
        W(x) &\mid x^{2a + 2b} + 1 ,
    \end{align*}
    which promptly implies
    \begin{alignat*}{2}
        && W(x) &\mid (1+x^{12b}-x^{a+7b}) (x^{2a + b} + x^a + x^{2b})\\
        && & \hspace{0.9cm} + (-x^b + x^{a+2b} - x^{a+11b} + x^{2a+6b})(x^{a + 2b} + x^b + x^{2a})\\
        && & \hspace{0.9cm} + (-x^a + x^{a + 9b} - x^{2a + 4b})(x^{2a + 2b} + 1)\\
        \implies \quad && W(x) &\mid x^{14b} ,
    \end{alignat*}
    which is again not possible.
\end{proof}

\begin{lemma}\label{complicated_lemma_2}
    If a given polynomial $W(x) \in \mathbb{Q}[x]$ has at least two nonzero terms and all of its nonzero terms have powers divisible by four, then $W(x) \mid Q_{a, b}(x)$ and $W(x) \mid R_{a, b}(x)$ both imply $2 \mid a, b$.
\end{lemma}
\begin{proof}
    If $2 \mid a, b$, then the lemma statement holds and there is nothing left to prove. Given the fact that $Q_{a, b}(x) = Q_{b, a}(x)$ and $R_{a, b}(x) = R_{b, a}(x)$, we may assume, without loss of generality, that $a$ is odd. Depending on the value of $b \bmod 4$, we can divide the problem into four corresponding cases and solve each of them separately.

    \begin{table}[h!t]
    {\footnotesize
    \begin{center}
    \begin{tabular}{rcccc}
    \toprule & $b \equiv_4 2$ & $b \equiv_4 0$ & $b \equiv_4 a$ & $b \equiv_4 a + 2$ \\
    \midrule
    $2a+b \equiv_4$ & $0$ & $2$ & $a+2$ & $a$\\
    $a+2b \equiv_4$ & $a$ & $a$ & $a+2$ & $a+2$\\
    $a \equiv_4$ & $a$ & $a$ & $a$ & $a$\\
    $b \equiv_4$ & $2$ & $0$ & $a$ & $a+2$\\
    $2a+2b \equiv_4$ & $2$ & $2$ & $0$ & $0$\\
    $2a \equiv_4$ & $2$ & $2$ & $2$ & $2$\\
    $2b \equiv_4$ & $0$ & $0$ & $2$ & $2$\\
    $0 \equiv_4$ & $0$ & $0$ & $0$ & $0$\\
    \bottomrule
    \end{tabular}
    \end{center}
    \caption{The powers of the $Q_{a, b}(x)$ and $R_{a, b}(x)$ nonzero terms modulo four, provided $a$ is odd.}
    \label{complicated_remainders_2}
    }
    \end{table}

    \bigskip\noindent
    \textbf{Case 1:} $b \equiv_4 2$.\quad
    In this case, Table \ref{complicated_remainders_2} and Lemma \ref{modulo_div_lemma} tell us that $W(x) \mid Q_{a, b}(x)$ implies
    \begin{align}
        \label{aux_26} W(x) &\mid x^b - x^{2a + 2b} - x^{2a},\\
        \label{aux_27} W(x) &\mid x^{a + 2b} + x^a.
    \end{align}
    However, by combining Equations (\ref{aux_26}) and (\ref{aux_27}), we get
    \begin{alignat*}{2}
        && W(x) &\mid (x^b - x^{2a + 2b} - x^{2a}) + x^a \, (x^{a + 2b} + x^a)\\
        \implies \quad && W(x) &\mid x^b ,
    \end{alignat*}
    which is not possible. Similarly, $W(x) \mid R_{a, b}(x)$ implies
    \begin{align}
        \label{aux_28} W(x) &\mid x^b + x^{2a + 2b} + x^{2a},\\
        \label{aux_29} W(x) &\mid x^{a + 2b} + x^a,
    \end{align}
    which can also be shown to be impossible by combining Equation (\ref{aux_28}) and (\ref{aux_29})
    \begin{alignat*}{2}
        && W(x) &\mid (x^b + x^{2a + 2b} + x^{2a}) - x^a \, (x^{a + 2b} + x^a)\\
        \implies \quad && W(x) &\mid x^b .
    \end{alignat*}

    \bigskip\noindent
    \textbf{Case 2:} $b \equiv_4 0$.\quad
    Here, Table \ref{complicated_remainders_2} and Lemma \ref{modulo_div_lemma} dictate that $W(x) \mid Q_{a, b}(x)$ implies
    \begin{align*}
        W(x) &\mid x^{a + 2b} + x^a,\\
        W(x) &\mid x^{2a + b} - x^{2a + 2b} - x^{2a},
    \end{align*}
    while $W(x) \mid R_{a, b}(x)$ implies
    \begin{align*}
        W(x) &\mid x^{a + 2b} + x^a,\\
        W(x) &\mid x^{2a + b} + x^{2a + 2b} + x^{2a} .
    \end{align*}
    This case can now be resolved in an entirely analogous manner as done so in case~3 from the proof of Lemma \ref{complicated_lemma_p}.

    \bigskip\noindent
    \textbf{Case 3:} Case $b \equiv_4 a$.\quad
    Here, Table \ref{complicated_remainders_2} and Lemma \ref{modulo_div_lemma} dictate that $W(x) \mid Q_{a, b}(x)$ and $W(x) \mid R_{a, b}(x)$ would both have to imply
    \begin{align}
        \label{aux_30} W(x) &\mid x^a + x^b,\\
        \label{aux_31} W(x) &\mid x^{2a} + x^{2b}.
    \end{align}
    By combining Equations (\ref{aux_30}) and (\ref{aux_31}) we could obtain
    \begin{alignat*}{2}
        && W(x) &\mid (x^{2a} + x^{2b}) + (x^a - x^b)(x^a + x^b)\\
        \implies \quad && W(x) &\mid 2x^{2a} ,
    \end{alignat*}
    which is not possible, as desired.

    \bigskip\noindent
    \textbf{Case 4:} $b \equiv_4 a + 2$.\quad
    In this situation, from Table \ref{complicated_remainders_2} and Lemma \ref{modulo_div_lemma} we can easily see that both $W(x) \mid Q_{a, b}(x)$ and $W(x) \mid R_{a, b}(x)$ would surely imply
    \begin{align}
        \label{aux_32} W(x) &\mid x^{2a+b} + x^a,\\
        \label{aux_33} W(x) &\mid x^{2a+2b} + 1.
    \end{align}
    However, if we combined Equations (\ref{aux_32}) and (\ref{aux_33}), we would get
    \begin{alignat*}{2}
        && W(x) &\mid (x^{a + b} - 1)(x^{2a + b} + x^a) + x^a \, (x^{2a + 2b} + 1)\\
        \implies \quad && W(x) &\mid 2x^{3a + 2b} ,
    \end{alignat*}
    which is impossible once again.
\end{proof}

It is worth pointing out that the polynomials $Q_{a, b}(x)$ and $R_{a, b}(x)$ have at most eight nonzero terms. For this reason, it becomes convenient to implement Theorem~\ref{filaseta} while inspecting their divisibility by cyclotomic polynomials. For example, if $\Phi_f(x) \mid Q_{a, b}(x)$ for some $f \in \mathbb{N}$, we may cancel out each prime factor of $f$ which is greater than seven in order to obtain another integer $f' \in \mathbb{N}$ which satisfies $\Phi_{f'}(x) \mid Q_{a, b}(x)$. The same conclusion can be made regarding the $R_{a, b}(x)$ polynomials. The proof of Theorem \ref{type_a_th} will heavily rely on this observation. Before we proceed with the main part of the proof, we shall need one more auxiliary lemma regarding the divisibility of $Q_{a, b}(x)$ and $R_{a, b}(x)$ by certain cyclotomic polynomials.

\begin{lemma}\label{div_p_cyclotomic}
    For each prime $p \ge 11$, $\Phi_p(x) \mid Q_{a, b}(x)$ and $\Phi_p(x) \mid R_{a, b}(x)$ both imply $p \mid a, b$.
\end{lemma}
\begin{proof}
    Let us define $Q_{a, b}^{\bmod p}(x)$ and $R_{a, b}^{\bmod p}(x)$ as the polynomials obtained from $Q_{a, b}(x)$ and $R_{a, b}(x)$, respectively, by replacing the power of each nonzero term by its remainder modulo $p$. It is straightforward to see that $\Phi_p(x) \mid Q_{a, b}(x)$ holds if and only if $\Phi_p(x) \mid Q_{a, b}^{\bmod p}(x)$ does. Also, $\Phi_p(x) \mid R_{a, b}(x)$ is surely equivalent to $\Phi_p(x) \mid R_{a, b}^{\bmod p}(x)$.
    
    We shall demonstrate the lemma statement only for the $Q_{a, b}(x)$ polynomial, given the fact that the proof regarding the $R_{a, b}(x)$ polynomial can be carried out in an entirely analogous manner. Suppose that $\Phi_p(x) \mid Q_{a, b}(x)$. Bearing in mind that
    \[
        \Phi_p(x) = \sum_{j = 0}^{p-1} x^j ,
    \]
    it is easy to see that $\deg \Phi_p(x) = p - 1$. However, since $\deg Q_{a, b}^{\bmod p}(x) \le p-1$, the divisibility $\Phi_p(x) \mid Q_{a, b}^{\bmod p}(x)$ directly implies that one of the following two possibilities must be true:
    \begin{itemize}
        \item $Q_{a, b}^{\bmod p}(x) = \beta \, \Phi_p(x)$ for some $\beta \in \mathbb{Q} \setminus \{ 0 \}$;
        \item $Q_{a, b}^{\bmod p}(x) \equiv 0$.
    \end{itemize}
    In the first scenario, $Q_{a, b}^{\bmod p}(x)$ would need to have exactly $p \ge 11$ nonzero terms, which is obviously not possible. Thus, $Q_{a, b}^{\bmod p}(x) \equiv 0$ surely holds. It is now convenient to perform a remainder modulo $p$ analysis identical to the one done throughout the proof of Lemma \ref{complicated_lemma_p}. Bearing in mind all the cases disclosed in the aforementioned lemma, we may conclude that at least one of the following four statements must be true:
    \begin{itemize}
        \item $p \mid a, b$;
        \item there exists an element from the set $\{2a + b, a + 2b, a, b, 2a + 2b, 2a, 2b, 0 \}$ whose remainder modulo $p$ is unique within that set;
        \item the numbers $\{2a + b, a + 2b, a, b, 2a + 2b, 2a, 2b, 0 \}$ can be partitioned as $\{ \{2a + b, b \}, \{a + 2b, 2a + 2b, 2b \}, \{ a, 2a, 0 \} \}$ according to their remainder modulo $p$.
        \item the numbers $\{2a + b, a + 2b, a, b, 2a + 2b, 2a, 2b, 0 \}$ can be partitioned as $\{ \{a + 2b, a \}, \{2a + b, 2a + 2b, 2a \}, \{ b, 2b, 0 \} \}$ according to their remainder modulo $p$.
    \end{itemize}
    If $p \mid a, b$, then there is nothing left to discuss, since the lemma statement directly holds. If the set $\{2a + b, a + 2b, a, b, 2a + 2b, 2a, 2b, 0 \}$ contains an element whose remainder modulo $p$ is unique within the set, it is easy to see that $Q_{a, b}^{\bmod p}(x) \equiv 0$ leads to a contradiction. The last two scenarios can be resolved in an analogous manner, so we will only deal with the third.

    Suppose that the elements of the set $\{2a + b, a + 2b, a, b, 2a + 2b, 2a, 2b, 0 \}$ can be partitioned as $\{ \{2a + b, b \}, \{a + 2b, 2a + 2b, 2b \}, \{ a, 2a, 0 \} \}$ according to their remainder modulo $p$. From here, the condition $Q_{a, b}^{\bmod p}(x) \equiv 0$ immediately implies
    \[
        x^{a \bmod p} - x^{2a \bmod p} - 1 \equiv 0,
    \]
    which is clearly impossible.
\end{proof}

We are now finally able to put all the pieces of the puzzle together and finalize the proof of Theorem \ref{type_a_th}.

\bigskip\noindent
\emph{Proof of Theorem \ref{type_a_th}}.\quad
Let $G$ be an arbitrarily chosen graph which is representable as $T_4(n, a, b)$. According to Lemma \ref{type_a_easy_direction}, if $G$ does not satisfy the four conditions given in the theorem, then this graph is certainly not a nut graph. Thus, in order to complete the proof, it is sufficent to suppose that $G$ is not a nut graph, then prove that it fails to satisfy at least one of the four stated conditions. However, if we suppose that the graph $G$ is not a nut graph, then Lemma \ref{type_a_cyclotomic_lemma} states that there must exist an $f \in \mathbb{N}, \, f \ge 2$ such that $f \mid \frac{n}{2}$ and $\Phi_f(x) \mid Q_{a, b}(x)$, or an $f \in \mathbb{N}$ such that $f \mid n$, $2 \nmid \frac{n}{f}$ and $\Phi_f(x) \mid R_{a, b}(x)$. It now becomes convenient to split the problem into two cases depending on whether $Q_{a, b}(x)$ or $R_{a, b}(x)$ is divisible by the corresponding cyclotomic polynomial.

\bigskip\noindent
\textbf{Case 1}: $\Phi_f(x) \mid Q_{a, b}(x)$.\quad
If $p^2 \mid f$ holds for any prime number $p \in \mathbb{N}$, we then have that $\Phi_f(x) = \Phi_{f/p}(x^p)$, which implies that all the nonzero terms of $\Phi_f(x)$ possess powers divisible by $p$. Thus, if $p^2 \mid f$ for any prime $p \ge 7$ or $p = 3$, we can implement Lemma \ref{complicated_lemma_p} or Lemma \ref{complicated_lemma_3} in order to obtain that $p \mid a, b$. However, it is now straightforward to see that $p \mid a, b, \frac{n}{2}$, hence $\gcd(\frac{n}{2}, a, b) \neq 1$, which means that condition \textbf{(i)} from Theorem \ref{type_a_th} is not satisfied.

If we have that $5^2 \mid f$, we may apply Lemma \ref{complicated_lemma_5} to obtain that $5 \mid a, b$ or $5 \nmid a, b, a-b, a+b$. If $5 \mid a, b$, it is then clear that condition \textbf{(i)} from Theorem~\ref{type_a_th} fails to hold. If $5 \nmid a, b, a-b, a+b$, it is enough to notice that $10 \mid n$ to conclude that condition \textbf{(iv)} from Theorem \ref{type_a_th} is not satisfied. Also, if $8 \mid f$, it is then convenient to use Lemma \ref{complicated_lemma_2} to reach $2 \mid a, b$. However, this further gives that condition \textbf{(i)} from Theorem~\ref{type_a_th} fails to hold once again.

Taking everything into consideration, we may assume that $\Phi_f(x) \mid Q_{a, b}(x)$ where $f \in \mathbb{N}$ is such that:
\begin{itemize}
    \item $f \ge 2$ and $f \mid \frac{n}{2}$;
    \item $8 \nmid f$ and $p^2 \nmid f$ for every prime $p \ge 3$.
\end{itemize}
We now divide the problem into two subcases depending on whether $f$ contains a prime factor from the set $\{2, 3, 5, 7 \}$.

\medskip\noindent
\textbf{Subcase 1.1:} $2, 3, 5, 7 \nmid f$.\quad
In this subcase, it is obvious that $f$ must be representable as a product of one or more distinct prime numbers greater than seven. By applying Theorem \ref{filaseta}, we may cancel out all the prime factors of $f$, one by one, until there is exactly one left. This allows us to conclude that $\Phi_p(x) \mid Q_{a, b}(x)$ for some prime $p \ge 11$ such that $p \mid f$, hence $p \mid \frac{n}{2}$. By implementing Lemma \ref{div_p_cyclotomic}, it becomes easy to see that $p \mid a, b$. However, this means that condition \textbf{(i)} from Theorem \ref{type_a_th} does not hold.

\medskip\noindent
\textbf{Subcase 1.2:} $\lnot (2, 3, 5, 7 \nmid f)$.\quad
Here, we have that $f$ contains at least one prime factor not greater than seven. We are now able to use Theorem \ref{filaseta} in order to cancel out all the prime factors of $f$ greater than seven. Furthermore, if the obtained integer is simultaneously divisible by five and seven, we can implement Theorem \ref{filaseta} once more and cancel out one of them, given the fact that $(5-2) + (7-2) > (8-2)$. Bearing everything in mind, we get that $\Phi_{f'}(x) \mid Q_{a, b}(x)$ must hold for some integer $f' \in \mathbb{N}$ such that:
\begin{itemize}
    \item $f' \ge 2$, $f' \mid \frac{n}{2}$ and $f'$ contains no prime factor greater than seven;
    \item $8 \nmid f'$ and $p^2 \nmid f'$ for each prime $p \in \{3, 5, 7\}$;
    \item $f'$ is not divisible by both five and seven.
\end{itemize}
It is not difficult to check that such an integer $f'$ would necessarily have to belong to the set
\[
    \{ 2, 3, 4, 5, 6, 7, 10, 12, 14, 15, 20, 21, 28, 30, 42, 60, 84 \} .
\]

Furthermore, let us define $Q_{a, b}^{\bmod f'}(x)$ as the polynomial obtained from $Q_{a, b}(x)$ by replacing the power of each nonzero term by its remainder modulo $f'$. It is clear that $\Phi_{f'}(x) \mid Q_{a, b}(x)$ is equivalent to $\Phi_{f'}(x) \mid Q_{a, b}^{\bmod f'}(x)$. We now observe that for a fixed $f' \in \mathbb{N}$ one can determine a finite set
\[
    \Psi \subseteq \{ 0, 1, \ldots, n-1 \} \times \{0, 1, \ldots, n-1 \},
\]
such that $\Phi_{f'}(x) \mid Q_{a, b}^{\bmod f'}(x)$ holds if and only if $(a \bmod p, b \bmod p) \in \Psi$. Given the fact that there are only $17$ possible values for $f'$, it is quite convenient to use a computer in order to determine all the modular conditions under which $\Phi_{f'}(x) \mid Q_{a, b}(x)$ is true for each feasible $f'$. The corresponding computational results are given in Appendix \ref{q_inspection}. By examining the said results, it is not difficult to establish that at least one of the following statements has to be true:
\begin{itemize}
    \item there exists a prime number $p \in \mathbb{N}$ such that $p \mid a, b, f'$;
    \item $f' = 4$ and $a$ and $b$ are both odd;
    \item $f' = 5$ and $5 \nmid a, b, a+b, a-b$.
\end{itemize}
If $p \mid a, b, f'$, we then clearly have $\gcd(\frac{n}{2}, a, b) \neq 1$, which implies that condition~\textbf{(i)} given in Theorem \ref{type_a_th} does not hold. Similarly, if $f' = 4$ and $2 \nmid a, b$, it can be immediately seen that $4 \mid \frac{n}{2}$, hence condition \textbf{(iii)} is not satisfied. Finally, if $f' = 5$ and $5 \nmid a, b, a+b, a-b$, then $5 \mid \frac{n}{2}$, which means that condition \textbf{(iv)} from Theorem \ref{type_a_th} does not hold, as desired.

\bigskip\noindent
\textbf{Case 2:} $\Phi_f(x) \mid R_{a, b}(x)$.\quad
In this case, we can apply the same initial discussion once again in order to show that we may asume that $\Phi_f(x) \mid R_{a, b}(x)$ holds for some integer $f \in \mathbb{N}$ such that:
\begin{itemize}
    \item $f \mid n$ and $2 \nmid \frac{n}{f}$;
    \item $8 \nmid f$ and $p^2 \nmid f$ for every prime $p \ge 3$.
\end{itemize}
It is now convenient to divide the problem into two subcases in the same manner as it was done in the previous case.

\medskip\noindent
\textbf{Subcase 2.1:} $2, 3, 5, 7 \nmid f$.\quad
This subcase can be resolved in an entirely analogous manner as subcase 1.1. For this reason, we choose to omit the proof details.

\medskip\noindent
\textbf{Subcase 2.2:} $\lnot (2, 3, 5, 7 \nmid f)$.\quad
In this subcase, we can implement Theorem \ref{filaseta} in an identical manner as done so in subcase 1.2. This allows us to reach that $\Phi_{f'}(x) \mid R_{a, b}(x)$ must be true for some integer $f' \in \mathbb{N}$ such that:
\begin{itemize}
    \item $f' \mid n$, $2 \nmid \frac{n}{f'}$ and $f'$ contains no prime factor greater than seven;
    \item $8 \nmid f'$ and $p^2 \nmid f'$ for each prime $p \in \{3, 5, 7\}$;
    \item $f'$ is not divisible by both five and seven.
\end{itemize}
Taking into account that $f' \neq 1$, this means that $f'$ once again certainly belongs to the set
\[
    \{ 2, 3, 4, 5, 6, 7, 10, 12, 14, 15, 20, 21, 28, 30, 42, 60, 84 \} .
\]

By defining $R_{a, b}^{\bmod p}(x)$ in an analogous manner as $Q_{a, b}^{\bmod p}(x)$ was defined in case~1, it becomes possible to inspect the precise modular conditions that $a$ and $b$ have to satisfy in order for $\Phi_{f'} \mid R_{a, b}(x)$ to be true for a given value of $f'$. Of course, the aforementioned examination can be easily performed via computer. The corresponding computational results are disclosed in Appendix \ref{r_inspection} and by analyzing the obtained results it is possible to conclude that $f'$ is surely even. Moreover, at least one of the following statements is certainly true:
\begin{itemize}
    \item there exists a prime number $p \in \mathbb{N}$ such that $p \mid a, b, \frac{f'}{2}$;
    \item $f' = 2$ and $a$ and $b$ are both odd;
    \item $f' = 4$ and $a$ and $b$ are both odd;
    \item $f' = 10$ and $5 \nmid a, b, a-b, a+b$.
\end{itemize}
If there exists a prime $p$ such that $p \mid a, b, \frac{f'}{2}$, it is evident that $\gcd(\frac{n}{2}, a, b) \neq 1$, hence condition \textbf{(i)} from Theorem \ref{type_a_th} is not satisfied. Furthermore, if $f' = 2$ and $2 \nmid a, b$, then it is easy to see that $n \equiv_4 2$, which means that condition \textbf{(ii)} fails to hold. Similarly, if $f' = 4$ and $2 \nmid a, b$, then we get $4 \mid n$, which implies that condition \textbf{(iii)} is not satisfied. Finally, if $f' = 10$ and $5 \nmid a, b, a-b, a+b$, it is straightforward to deduce that condition \textbf{(iv)} from Theorem \ref{type_a_th} fails to hold. \hfill\qed

\begin{remark}
    In an analogous manner as done so in Section \ref{sc_type_b}, it is possible to demonstrate that $T_4(n, a, b)$ is always isomorphic to
    \[
        T_4\left( n, \min(ta \bmod n, n - ta \bmod n), \min(tb \bmod n, n - tb \bmod n) \right)
    \]
    whenever $\gcd(n, t) = 1$. From here, it follows that if $4 \mid n$, then a tricirculant graph of type 4 is a nut graph if and only if it is representable as $T_4(n, a, b)$ where $1 \le a, b < \frac{n}{2}$, $\gcd(a, b) = 1$, $a \not\equiv_2 b$, and $5 \mid ab(a-b)(a+b)$ provided $10 \mid n$.
\end{remark}

\section*{Acknowledgements}

I.\ Damnjanović is grateful to the University of Primorska, University of Ljub\-ljana, Inštitut za matematiko, fiziko in mehaniko and Diffine LLC for all the support given throughout the duration of the research. T.\ Pisanski gratefully acknowledges the support provided by the Mathematical Institute of the Serbian Academy of Sciences and Arts. This work is supported in part by the Slovenian Research Agency (research program P1-0294 and research projects N1-0140, J1-2481, \linebreak J1-3002 and J1-4351).

\section*{Conflict of interest}

The authors declare that they have no conflict of interest.

\appendix

\section{Inspection for \texorpdfstring{$\Phi_f(x) \mid Q_{a, b}(x)$}{Phi(f)(x) | Q(a, b)(x)}}\label{q_inspection}

In this appendix section, we disclose the computational results that describe all the possible modular conditions that $a, b \in \mathbb{N}$ have to adhere to in order for $\Phi_f(x) \mid Q_{a, b}(x)$ to be true, for certain values of $f \in \mathbb{N}$. The given results can be generated, for example, by using the following Wolfram Mathematica command:

\begin{lstlisting}[language = Mathematica, frame = trBL, escapeinside={(*@}{@*)}, aboveskip=10pt, belowskip=10pt, numbers=left, rulecolor=\color{black}, morekeywords={CoefficientRules}]
MatrixForm[
 Table[{f, 
   MatrixForm[
    Select[Flatten[Table[{a, b}, {a, 0, f - 1}, {b, 0, f - 1}], 1], 
     Length[CoefficientRules[
         PolynomialRemainder[
          x^(2 #[[1]] + #[[2]]) + x^(#[[1]] + 2 #[[2]]) + x^#[[1]] + 
           x^#[[2]] - x^(2 #[[1]] + 2 #[[2]]) - x^(2 #[[1]]) - 
           x^(2 #[[2]]) - 1, Cyclotomic[f, x], x]]] == 0 &]]}, {f, {2,
     3, 4, 5, 6, 7, 10, 12, 14, 15, 20, 21, 28, 30, 42, 60, 84}}]]
\end{lstlisting}

\noindent\begin{center}
\noindent\resizebox{1.00\columnwidth}{!}{
{\scriptsize
\begin{tabular}[t]{l|ll}
\toprule $f$ & $a \bmod f$ & $b \bmod f$\\
\midrule
$2$ & $0$ & $0$\\
\midrule
$3$ & $0$ & $0$\\
\midrule
$4$ & $0$ & $0$\\
& $1$ & $1$ \\
& $1$ & $3$ \\
& $3$ & $1$ \\
& $3$ & $3$ \\
\midrule
$5$ & $0$ & $0$ \\
& $1$ & $2$ \\
& $1$ & $3$ \\
& $2$ & $1$ \\
& $2$ & $4$ \\
& $3$ & $1$ \\
& $3$ & $4$ \\
& $4$ & $2$ \\
& $4$ & $3$ \\
\midrule
$6$ & $0$ & $0$ \\
\midrule
$7$ & $0$ & $0$ \\
\midrule
$10$ & $0$ & $0$ \\
& $2$ & $4$ \\
& $2$ & $6$ \\
& $4$ & $2$ \\
& $4$ & $8$ \\
& $6$ & $2$ \\
& $6$ & $8$ \\
& $8$ & $4$ \\
& $8$ & $6$ \\
\midrule
$12$ & $0$ & $0$ \\
& $3$ & $3$ \\
& $3$ & $9$ \\
& $9$ & $3$ \\
& $9$ & $9$ \\
\bottomrule
\end{tabular}}
\hspace{0.8cm} {\scriptsize
\begin{tabular}[t]{l|ll}
\toprule $f$ & $a \bmod f$ & $b \bmod f$\\
\midrule
$14$ & $0$ & $0$ \\
\midrule
$15$ & $0$ & $0$ \\
& $3$ & $6$ \\
& $3$ & $9$ \\
& $6$ & $3$ \\
& $6$ & $12$ \\
& $9$ & $3$ \\
& $9$ & $12$ \\
& $12$ & $6$ \\
& $12$ & $9$ \\
\midrule
$20$ & $0$ & $0$ \\
& $4$ & $8$ \\
& $4$ & $12$ \\
& $5$ & $5$ \\
& $5$ & $15$ \\
& $8$ & $4$ \\
& $8$ & $16$ \\
& $12$ & $4$ \\
& $12$ & $16$ \\
& $15$ & $5$ \\
& $15$ & $15$ \\
& $16$ & $8$ \\
& $16$ & $12$ \\
\midrule
$21$ & $0$ & $0$ \\
\midrule
$28$ & $0$ & $0$ \\
& $7$ & $7$ \\
& $7$ & $21$ \\
& $21$ & $7$ \\
& $21$ & $21$ \\
\bottomrule
\end{tabular}}
\hspace{0.8cm} {\scriptsize
\begin{tabular}[t]{l|ll}
\toprule $f$ & $a \bmod f$ & $b \bmod f$\\
\midrule
$30$ & $0$ & $0$ \\
& $6$ & $12$ \\
& $6$ & $18$ \\
& $12$ & $6$ \\
& $12$ & $24$ \\
& $18$ & $6$ \\
& $18$ & $24$ \\
& $24$ & $12$ \\
& $24$ & $18$ \\
\midrule
$42$ & $0$ & $0$ \\
\midrule
$60$ & $0$ & $0$ \\
& $12$ & $24$ \\
& $12$ & $36$ \\
& $15$ & $15$ \\
& $15$ & $45$ \\
& $24$ & $12$ \\
& $24$ & $48$ \\
& $36$ & $12$ \\
& $36$ & $48$ \\
& $45$ & $15$ \\
& $45$ & $45$ \\
& $48$ & $24$ \\
& $48$ & $36$ \\
\midrule
$84$ & $0$ & $0$ \\
& $21$ & $21$ \\
& $21$ & $63$ \\
& $63$ & $21$ \\
& $63$ & $63$ \\
\bottomrule
\end{tabular}}
}
\end{center}

\section{Inspection for \texorpdfstring{$\Phi_f(x) \mid R_{a, b}(x)$}{Phi(f)(x) | R(a, b)(x)}}\label{r_inspection}

In this appendix section, we give the computational results that describe all the possible modular conditions that $a, b \in \mathbb{N}$ have to satisfy in order for $\Phi_f(x) \mid R_{a, b}(x)$ to hold, for concrete values of $f \in \mathbb{N}$. The said results can be quickly obtained, for example, by using the following Wolfram Mathematica command:

\pagebreak
\begin{lstlisting}[language = Mathematica, frame = trBL, escapeinside={(*@}{@*)}, aboveskip=10pt, belowskip=10pt, numbers=left, rulecolor=\color{black}, morekeywords={CoefficientRules}]
MatrixForm[
 Table[{f, 
   MatrixForm[
    Select[Flatten[Table[{a, b}, {a, 0, f - 1}, {b, 0, f - 1}], 1], 
     Length[CoefficientRules[
         PolynomialRemainder[
          x^(2 #[[1]] + #[[2]]) + x^(#[[1]] + 2 #[[2]]) + x^#[[1]] + 
           x^#[[2]] + x^(2 #[[1]] + 2 #[[2]]) + x^(2 #[[1]]) + 
           x^(2 #[[2]]) + 1, Cyclotomic[f, x], x]]] == 0 &]]}, {f, {2,
     3, 4, 5, 6, 7, 10, 12, 14, 15, 20, 21, 28, 30, 42, 60, 84}}]]
\end{lstlisting}

\noindent\begin{center}
\noindent\resizebox{1.00\columnwidth}{!}{
{\scriptsize
\begin{tabular}[t]{l|ll}
\toprule $f$ & $a \bmod f$ & $b \bmod f$ \\
\midrule
$2$ & $1$ & $1$\\
\midrule
$3$ & &\\
\midrule
$4$ & $1$ & $1$ \\
& $1$ & $3$ \\
& $2$ & $2$ \\
& $3$ & $1$ \\
& $3$ & $3$ \\
\midrule
$5$ & & \\
\midrule
$6$ & $3$ & $3$ \\
\midrule
$7$ & & \\
\midrule
$10$ & $1$ & $3$ \\
& $1$ & $7$ \\
& $3$ & $1$ \\
& $3$ & $9$ \\
& $5$ & $5$ \\
& $7$ & $1$ \\
& $7$ & $9$ \\
& $9$ & $3$ \\
& $9$ & $7$ \\
\midrule
$12$ & $3$ & $3$ \\
& $3$ & $9$ \\
& $6$ & $6$ \\
& $9$ & $3$ \\
& $9$ & $9$ \\
\bottomrule
\end{tabular}}
\hspace{0.8cm} {\scriptsize
\begin{tabular}[t]{l|ll}
\toprule $f$ & $a \bmod f$ & $b \bmod f$ \\
\midrule
$14$ & $7$ & $7$ \\
\midrule
$15$ & & \\
\midrule
$20$ & $2$ & $6$ \\
& $2$ & $14$ \\
& $5$ & $5$ \\
& $5$ & $15$ \\
& $6$ & $2$ \\
& $6$ & $18$ \\
& $10$ & $10$ \\
& $14$ & $2$ \\
& $14$ & $18$ \\
& $15$ & $5$ \\
& $15$ & $15$ \\
& $18$ & $6$ \\
& $18$ & $14$ \\
\midrule
$21$ & & \\
\midrule
$28$ & $7$ & $7$ \\
& $7$ & $21$ \\
& $14$ & $14$ \\
& $21$ & $7$ \\
& $21$ & $21$ \\
\bottomrule
\end{tabular}}
\hspace{0.8cm} {\scriptsize
\begin{tabular}[t]{l|ll}
\toprule $f$ & $a \bmod f$ & $b \bmod f$ \\
\midrule
$30$ & $3$ & $9$ \\
& $3$ & $21$ \\
& $9$ & $3$ \\
& $9$ & $27$ \\
& $15$ & $15$ \\
& $21$ & $3$ \\
& $21$ & $27$ \\
& $27$ & $9$ \\
& $27$ & $21$ \\
\midrule
$42$ & $21$ & $21$ \\
\midrule
$60$ & $6$ & $18$ \\
& $6$ & $42$ \\
& $15$ & $15$ \\
& $15$ & $45$ \\
& $18$ & $6$ \\
& $18$ & $54$ \\
& $30$ & $30$ \\
& $42$ & $6$ \\
& $42$ & $54$ \\
& $45$ & $15$ \\
& $45$ & $45$ \\
& $54$ & $18$ \\
& $54$ & $42$ \\
\midrule
$84$ & $21$ & $21$ \\
& $21$ & $63$ \\
& $42$ & $42$ \\
& $63$ & $21$ \\
& $63$ & $63$ \\
\bottomrule
\end{tabular}}
}
\end{center}


\begin{thebibliography}{xx}

\bibitem{BaKnSk2022}
N.\ Bašić, M.\ Knor, R.\ Škrekovski, On $12$-regular nut graphs, {\em Art Discrete Appl.\ Math.\/} {\bf 5(2)} (2022), \#P2.01, \doi{10.26493/2590-9770.1403.1b1}.

\bibitem{CoFoGo2018} K.\ Coolsaet, P.W.\ Fowler, J.\ Goedgebeur, Generation and properties of nut graphs, {\em MATCH Commun.\ Math.\ Comput.\ Chem.\/} {\bf 80} (2018), 423--444, URL: \href{https://match.pmf.kg.ac.rs/electronic_versions/Match80/n2/match80n2_423-444.pdf}{\texttt{match80n2\_423-444.pdf}}.

\bibitem{Cvetkovic1995}
D.M.\ Cvetković, M.\ Doob, H.\ Sachs, \emph{Spectra of graphs: Theory and applications}, Leipzig: J.\ A.\ Barth Verlag, 1995.

\bibitem{Damnjanovic2022_A} I.\ Damnjanović, Two families of circulant nut graphs, {\em Filomat} {\bf 37(24)} (2023), 8331--8360, \doi{10.2298/FIL2324331D}.

\bibitem{Damnjanovic2022_B} I.\ Damnjanović, Complete resolution of the circulant nut graph order--degree existence problem, {\em Ars Math. Contemp.\/} (2023), DOI: \href{https://doi.org/10.26493/1855-3974.3009.6df}{\texttt{10.26493/\linebreak1855-3974.3009.6df}}.

\bibitem{Damnjanovic2022_C} I.\ Damnjanović, On the nullities of quartic circulant graphs and their extremal null spaces, \arxiv{2212.12959}, 2022.

\bibitem{DaSt2022}
I.\ Damnjanović, D.\ Stevanović, On circulant nut graphs, {\em Linear Algebra Appl.\/} {\bf 633} (2022), 127--151, \doi{10.1016/j.laa.2021.10.006}.

\bibitem{FiSc2003} M.\ Filaseta, A.\ Schinzel, On testing the divisibility of lacunary polynomials by cyclotomic polynomials, {\em Math.\ Comput.\/} {\bf 73(246)} (2003), 957--965, URL:  \href{https://www.jstor.org/stable/4099813}{\texttt{https://www.jstor.org/stable/4099813}}.

\bibitem{FoGaGoPiSc2020} P.W.\ Fowler, J.B.\ Gauci, J.\ Goedgebeur, T.\ Pisanski, I.\ Sciriha, Existence of regular nut graphs for degree at most $11$, {\em Discuss.\ Math.\ Graph Theory} {\bf 40} (2020), 533--557, \doi{10.7151/dmgt.2283}.

\bibitem{FoPiToBoSc2014} P.W.\ Fowler, B.T.\ Pickup, T.Z.\ Todorova, M.\ Borg, I.\ Sciriha, Omni-conducting and omni-insulating molecules, {\em J.\ Chem.\ Phys.\/} {\bf 140(5)} (2014), 054115, \doi{10.1063/1.4863559}.

\bibitem{FoPiBa2021} P.W.\ Fowler, T.\ Pisanski, N.\ Bašić, Charting the space of chemical nut graphs, {\em MATCH Commun.\ Math.\ Comput.\ Chem.\/} {\bf 86(3)} (2021), 519--538, URL: \href{https://match.pmf.kg.ac.rs/electronic_versions/Match86/n3/match86n3_519-538.pdf}{\texttt{match86n3\_519-538.pdf}}.

\bibitem{GaPiSc2023} J.B.\ Gauci, T.\ Pisanski, I.\ Sciriha, Existence of regular nut graphs and the Fowler construction, {\em Appl.\ Anal.\ Discrete Math.\/} (2023), DOI: \href{https://www.doi.org/10.2298/AADM190517028G}{\texttt{10.2298/\allowbreak AADM190517028G}}.

\bibitem{Gray} R.M.\ Gray, Toeplitz and circulant matrices: A review, {\em Found.\ Trends Commun.\ Inf.\/} Theory {\bf 2} (2006), 155--239, URL: \href{https://ee.stanford.edu/~gray/CIT006-journal.pdf}{\texttt{CIT006-journal.pdf}}.

\bibitem{KoKuMaWi2012} I.\ Kov\'{a}cs, K.\ Kutnar, D.\ Marušič, S.\ Wilson, Classification of cubic symmetric tricirculants, {\em Electron.\ J.\ Combin.\/} {\bf 19(2)} (2012), P24, \doi{10.37236/2371}.

\bibitem{Nagell1951} T.\ Nagell, \emph{Introduction to Number Theory}, Wiley, New York, 1951.

\bibitem{Pi2007} T.\ Pisanski, A classification of cubic bicirculants, {\em Discrete Math.\/} \textbf{307(3--5)} (2007), 567--578, \doi{10.1016/j.disc.2005.09.053}.

\bibitem{PiPoZi2022} T.\ Pisanski, P.\ Potočnik, A.\ Žitnik, Classification of quartic bicirculant nut graphs, manuscript.

\bibitem{PoTol2020} P.\ Potočnik, M. Toledo, Classification of cubic vertex-transitive tricirculants, {\em Ars Math.\ Contemp.\/} {\bf 18} (2020), 1--31, \doi{10.26493/1855-3974.1815.b52}.

\bibitem{Sciriha1997} I.\ Sciriha, On the coefficient of $\lambda$ in the characteristic polynomial of singular graphs, {\em Util.\ Math.\/} {\bf 52} (1997), 97--111.

\bibitem{Sciriha1998_A} I.\ Sciriha, On singular line graphs of trees, {\em Congr.\ Numerantium} {\bf 135} (1998), 73--91.

\bibitem{Sciriha1998_B} I.\ Sciriha, On the construction of graphs of nullity one, {\em Discrete Math.\/} {\bf 181(1--3)} (1998), 193--211, \doi{10.1016/S0012-365X(97)00036-8}.

\bibitem{Sciriha1999} I.\ Sciriha, The two classes of singular line graphs of trees, {\em Rend.\ Semin.\ Mat.\ Messina}, Ser.\ II {\bf 20(5)} (1999), 167--180.

\bibitem{Sciriha2007} I.\ Sciriha, A characterization of singular graphs, {\em Electron.\ J.\ Linear Algebra}, {\bf 16} (2007), 451--462, URL: \url{https://eudml.org/doc/129125}.

\bibitem{Sciriha2008} I.\ Sciriha, Coalesced and embedded nut graphs in singular graphs, {\em Ars Math.\ Contemp.\/} {\bf 1} (2008), 20--31, \doi{10.26493/1855-3974.20.7cc}.

\bibitem{ScFa2021} I.\ Sciriha, A.\ Farrugia, \emph{From nut graphs to molecular structure and conductivity}, Mathematical chemistry monographs, No.\ 23, University of Kragujevac, Kragujevac, 2021.

\bibitem{ScFo2007} I.\ Sciriha, P.W.\ Fowler, Nonbonding orbitals in fullerenes: Nuts and cores in singular polyhedral graphs, {\em J.\ Chem.\ Inf.\ Model.\/} {\bf  47(5)} (2007), 1763--1775, \doi{10.1021/ci700097j}.

\bibitem{ScFo2008} I.\ Sciriha, P.W.\ Fowler, On nut and core singular fullerenes, {\em Discrete Math.\/} {\bf 308(2--3)} (2008), 267--276, \doi{10.1016/j.disc.2006.11.040}.

\bibitem{ScGu1998} I.\ Sciriha, I.\ Gutman, Nut graphs: maximally extending cores, {\em Util.\ Math.\/} {\bf 54} (1998), 257--272.
		
\bibitem{Tattersall} J.\ Tattersall, 
\emph{Elementary number theory in nine chapters}, Cambridge University Press, New York, 1999. 

\bibitem{West} D.B.\ West, \emph{Introduction to Graph Theory}, second edition, Prentice Hall, 2000.

\bibitem{Cyclotomic} Cyclotomic polynomials, Encyclopedia of Mathematics, Springer Verlag GmbH, European Mathematical Society, URL: \href{https://encyclopediaofmath.org/index.php?title=Cyclotomic_polynomials}{\texttt{https://encyclopediaof\linebreak math.org/index.php?title=Cyclotomic\_polynomials}}.




%

\end{thebibliography}
\end{document}